\documentclass[11pt]{amsart}
\renewcommand{\bar}{\overline}
\usepackage{amsmath,amsthm,amscd,euscript}
\renewcommand{\bar}{\overline}
\def \r{\mathbb R}
\def \c{\mathbb C}
\def \q{\mathbb Q}
\def \z{\mathbb Z}

\DeclareMathOperator{\MOD}{mod}
\DeclareMathOperator{\isin}{lsin}
\DeclareMathOperator{\icos}{lcos}
\DeclareMathOperator{\itan}{ltan}
\DeclareMathOperator{\iarctan}{larctan}
\DeclareMathOperator{\il}{l\ell} \DeclareMathOperator{\is}{lS}
\DeclareMathOperator{\sgn}{sgn} \DeclareMathOperator{\lcm}{lcm}
\DeclareMathOperator{\sign}{sign}

\DeclareMathOperator{\pcong}{\hat{\cong}}

\newtheorem{theorem}{Theorem}[section]
\newtheorem{lemma}[theorem]{Lemma}
\newtheorem{statement}[theorem]{Statement}
\newtheorem{proposition}[theorem]{Proposition}
\newtheorem{corollary}[theorem]{Corollary}
\theoremstyle{remark}
\newtheorem{remark}[theorem]{Remark}
\theoremstyle{definition}
\newtheorem{definition}[theorem]{Definition}
\newtheorem{example}[theorem]{Example}
\newtheorem{problem}{Problem}

\title{Elementary Notions of Lattice Trigonometry.}
\author{Oleg~Karpenkov}
\date{22 March 2006}
\thanks{Partially supported
by NWO-RFBR 047.011.2004.026 (RFBR 05-02-89000-NWO\_a) grant, by
RFBR SS-1972.2003.1 grant, by RFBR 05-01-02805-CNRSL\_a grant,
and by RFBR grant 05-01-01012a.}


\email[Oleg Karpenkov]{karpenk@mccme.ru}

\address{CEREMADE - UMR 7534
Universit\'e Paris-Dauphine France -- 75775 Paris SEDEX 16}

\begin{document}
\input epsf

\begin{abstract}
In this paper we study properties of lattice trigonometric
functions of lattice angles in lattice geometry. We introduce the
definition of sums of lattice angles and establish a necessary
and sufficient condition for three angles to be the angles of
some lattice triangle in terms of lattice tangents. This
condition is a version of the Euclidean condition: three angles
are the angles of some triangle iff their sum equals $\pi$.
Further we find the necessary and sufficient condition for an
ordered $n$-tuple of angles to be the angles of some convex
lattice polygon. In conclusion we show applications to theory of
complex projective toric varieties, and a list of unsolved
problems and questions.
\end{abstract}

 \maketitle

\tableofcontents

\sloppy \normalsize

\section*{Introduction.}

\subsection{The goals of this paper and some background.}
Consider a two-dimensional oriented real vector space and fix
some full-rank lattice in it. A triangle or a polygon is said to
be {\it lattice} if all its vertices belong to the lattice. The
angles of any lattice triangle are said to be {\it lattice}.

In this paper we introduce and study {\it lattice trigonometric}
functions of lattice angles. The lattice trigonometric functions
are invariant under the action of the group of lattice-affine
transformations (i.e. affine transformations preserving the
lattice), like the ordinary trigonometric functions are invariant
under the action of the group of Euclidean length preserving
transformations of Euclidean space.

One of the initial goals of the present article is to make a
complete description of lattice triangles up to the
lattice-affine equivalence relation (see Theorem~\ref{sum}). The
classification problem of convex lattice polygons becomes now
classical. There is still no a good description of convex
polygons. It is only known that the number of such polygons with
lattice area bounded from above by $n$ growths exponentially in
$n$, while $n$ tends to infinity (see the works of
V.~Arnold~\cite{Arn5}, and of I.~B\'ar\'any and
A.~M.~Vershik~\cite{Bar2}).

We expand the geometric interpretation of ordinary continued
fractions to define lattice sums of lattice angles and to
establish relations on lattice tangents of lattice angles.
Further, we describe lattice triangles in terms of {\it lattice
sums} of lattice angles.

In present paper we also show a lattice version of the sine
formula and introduce a relation between the lattice tangents for
angles of lattice triangles  and the numbers of lattice points on
the edges of triangles (see Theorem~\ref{ang-length}). In
addition, we make first steps in study of non-lattice angles with
lattice vertices. We conclude the paper with applications to
toric varieties and some unsolved problems.

\vspace{2mm}

The study of lattice angles is an imprescriptible part of modern
lattice geometry. Invariants of lattice angles are used in the
study of lattice convex polygons and polytopes. Such polygons and
polytopes play the principal role in Klein's theory of
multidimensional continued fractions (see, for example, the works
of F.~Klein~\cite{Kle1}, V.~I.~Arnold~\cite{Arn2},
E.~Korkina~\cite{Kor2}, M.~Kontsevich and Yu.~Suhov~\cite{Kon},
G.~Lachaud~\cite{Lac2}, and the author~\cite{Kar1}).

Lattice polygons and polytopes of the lattice geometry are in the
limelight of complex projective toric varieties (see for more
information the works of V.~I.~Danilov~\cite{Dan},
G.~Ewald~\cite{Ewa}, T.~Oda~\cite{Oda}, and W.~Fulton~\cite{Ful}).
To illustrate, we deduce (in Appendix~\ref{toric}) from
Theorem~\ref{sum} the corresponding global relations on the toric
singularities for projective toric varieties associated to
integer-lattice triangles. We also show the following simple
fact: for any collection with multiplicities of
complex-two-dimensional toric algebraic singularities there
exists a complex-two-dimensional toric projective variety with
the given collection of toric singularities $($this result seems
to be classical, but it is missing in the literature$)$.

The studies of lattice angles and measures related to them were
started by A.~G.~Khovanskii, A.~Pukhlikov in~\cite{Kh1}
and~\cite{Kh2} in 1992. They introduced and investigated special
additive polynomial measure for the expanded notion of polytopes.
The relations between sum-formulas of lattice trigonometric
functions and lattice angles in Khovanskii-Pukhlikov sense are
unknown to the author.

\subsection{Some distinctions between lattice and Euclidean cases.} This paper is organized as follows.
Lattice trigonometric functions and Euclidean trigonometric
functions have much in common. For example, the values of lattice
tangents and Euclidean tangents coincide in a special natural
system of coordinates. Nevertheless, lattice geometry differs a
lot from Euclidean geometry.
We provide this with the following $4$ examples.\\
{\bf 1.} The angles $\angle ABC$ and $\angle CBA$ are always
congruent in Euclidean geometry,
but not necessary lattice-congruent in lattice geometry.\\
{\bf 2.} In Euclidean geometry for any $n\ge 3$ there exist a
regular polygon with $n$ vertices, and any two regular polygons
with the same number of vertices are homothetic to each other. In
lattice geometry there are only 6 non-homothetic regular lattice
polygons: two triangles (distinguished by lattice tangents of
angles), two quadrangles, and two octagons.
(See a more detailed description in~\cite{KarReg}.)\\
{\bf 3.} Consider three Euclidean criterions of triangle
congruence. Only the first criterion can be taken to the case of
lattice geometry. The others two are false in lattice
trigonometry.
(We refer to Appendix~\ref{criterio}.)\\
{\bf 4.} There exist two non-congruent right angles in lattice
geometry. (See Corollary~\ref{right_}.)

\subsection{Description of the paper.} This paper is organized as follows.

We start in Section~1 with some general notation of lattice
geometry and ordinary continued fractions. We define ordinary
lattice angles, and the functions lattice sine, tangent, and
cosine on the set of ordinary lattice angles, and lattice
arctangent for rationals greater than $1$. Further we indicate
their basic properties. We proceed with the geometrical
interpretation of lattice tangents in terms of ordinary continued
fractions. We also say a few words about relations of transpose
and adjacent lattice angles. In conclusion of Section~1 we study
of the basic properties of angles in lattice triangles.

In Section~2 we introduce the sum formula for the lattice
tangents of ordinary lattice angles of lattice triangles. The sum
formula is a lattice generalization of the following Euclidean
statement: three angles are the angles of some triangle iff their
sum equals $\pi$.

Further in Section~3 we introduce the notion of expanded lattice
angles and their normal forms and give the definition of sums of
expanded and ordinary lattice angles. For the definition of
expanded lattice angles we expand the notion of sails in the
sense of Klein: we define and study oriented broken lines on the
unit distance from lattice points.

In Section~4 we finally prove the first statement of the theorem
on sums of lattice tangents for ordinary lattice angles in
lattice triangles. In this section we also describe some
relations between continued fractions for lattice oriented broken
lines and the lattice tangents for the corresponding expanded
lattice angles. Further we give a necessary and sufficient
condition for an ordered $n$-tuple of angles to be the angles of
some convex lattice polygon.

In Section~5 we generalize the notions of ordinary and expanded
lattice angles and their sums to the case of angles with lattice
vertices but not necessary lattice rays. We find normal forms and
extend the definition of lattice sums for a certain special case
of such angles (we call them {\it irrational}).

We conclude this paper with three appendices. In Appendix~A we
describe applications to theory of complex projective toric
varieties mentioned above. Further in Appendix~B we give a list
of unsolved problems and questions. Finally in Appendix~C we
formulate criterions of lattice congruence for lattice triangles.
These criterions lead to the complete list of lattice triangles
with small lattice area (not greater than 10).

\vspace{2mm}

{\bf Acknowledgement.} The author is grateful to V.~I.~Arnold for
constant attention to this work, A.~G.~Khovanskii,
V.~M.~Kharlamov, J.-M.~Kantor, D.~Zvonkine, and D.~Panov for
useful remarks and discussions, and Universit\'e Paris-Dauphine
--- CEREMADE  for the hospitality and excellent working
conditions.

\section{Definitions and elementary properties of lattice trigonometric functions.}

We start the section with general definitions of lattice geometry
and notions of ordinary continued fractions. We define the
functions lattice sine, tangent, and cosine on the set of
ordinary lattice angles and formulate their basic properties.

Then we describe a geometric interpretation of lattice
trigonometric functions in terms of ordinary continued fractions
associated with boundaries of convex hulls for the sets of
lattice points contained in angles. We also say a few words about
relations of transpose and adjacent lattice angles.

Further we introduce the sine formula for the lattice angles of
lattice triangles. Finally, we show how to find the lattice
tangents of all angles and the lattice lengths of all edges of
any lattice triangle, if we know the lattice lengths of two edges
and the lattice tangent of the angle between them.

\subsection{Preliminary notions and definitions.}
By $\gcd(n_1,\ldots,n_k)$ and by $\lcm(n_1,\ldots,n_k)$ we denote
the greater common divisor and the less common multiple of the
nonzero integers $n_1, \ldots, n_k$ respectively. Suppose that
$a$, $b$ be arbitrary integers, and $c$ be an arbitrary positive
integer. We write that $a\equiv b (\MOD c)$ if the reminders of
$a$ and $b$ modulo $c$ coincide.

\subsubsection{Lattice notation.}
Here we define the main objects of lattice geometry, their lattice
characteristics, and the relation of lattice-congruence.

Consider a two-dimensional oriented real vector space and fix
some lattice in it. A straight line is said to be {\it lattice}
if it contains at least two distinct lattice points. A ray is
said to be {\it lattice} if it's vertex is lattice, and the
straight line containing the ray is lattice. An angle (i.e. the
union of two rays with the common vertex) is said to be {\it
ordinary lattice} if the rays defining it are lattice. A segment
is called {\it lattice} if its endpoints are lattice points.

By a {\it convex polygon} we mean a convex hulls of a finite
number of points that do not lie in a straight line. The minimal
set A straight line $\pi$ is said to be {\it supporting} for a
convex polygon $P$, if the intersections of $P$ and $\pi$ is not
empty, and the whole polygon $P$ is contained in one of the
closed half-planes bounded by $\pi$. An intersection of any
polygon $P$ with its supporting hyperplane is called a {\it
vertex} or an {\it edge} of the polygon if the dimension of
intersection is zero, or one respectively.

A triangle (or convex polygon) is said to be {\it lattice} if all
it's vertices are lattice points.

The affine transformation is called {\it lattice-affine} if it
preserves the set of all lattice points. Consider two arbitrary
(not necessary lattice in the above sense) sets. We say that
these two sets are {\it lattice-congruent} to each other if there
exist a lattice-affine transformation of $\r^2$ taking the first
set to the second.

A lattice triangle is said to be {\it simple} if the vectors
corresponding to its edges generate the lattice.
\begin{definition}
The {\it lattice length} of a lattice segment $AB$ is the ratio
between the Euclidean length of $AB$ and the length of the basic
lattice vector for the straight line containing this segment. We
denote the lattice length by $\il(AB)$.
\\
By the (non-oriented) {\it lattice area} of the convex polygon
$P$ we will call the ratio of the Euclidean area of the polygon
and the area of any lattice simple triangle, and denote it by
$\is(P)$.
\end{definition}

Any two rays (straight lines) are lattice-congruent to each other.
Two lattice segments are lattice-congruent iff they have equal
lattice lengths. The lattice area of the convex polygon is
well-defined and is proportional to the Euclidean area of the
polygon.

\subsubsection{Finite ordinary continued fractions.}
For any finite sequence $(a_0, a_1, \ldots ,a_n)$ where the
elements $a_1, \ldots , a_n$ are positive integers and $a_0$ is
an arbitrary integer we associate the following rational number
$q$:
$$
\begin{array}{ccc}
q&=&a_0+\frac{\displaystyle 1}{\displaystyle a_1+\frac{1}{
\begin{array}{cc}
\ddots&\vdots\\
& a_{n-1}+\frac{1}{a_n}\\
\end{array}
}}
\end{array}
.
$$
This representation of the rational $q$ is called an {\it
ordinary continued fraction} for $q$ and denoted by $[a_0, a_1,
\ldots ,a_n]$. (In the literature is also in use the following
notation: $[a_0; a_1, \ldots ,a_n]$.) An ordinary continued
fraction $[a_0, a_1, \ldots ,a_n]$ is said to be {\it odd} if
$n{+}1$ is odd, and {\it even} if $n{+}1$ is even.

Note that if $a_n\ne 1$ then $[a_0, a_1, \ldots ,a_n]=[a_0, a_1,
\ldots ,a_n-1,1]$.

Now we formulate a classical theorem of ordinary continued
fractions theory.

\begin{theorem}
For any rational there exist exactly one odd ordinary continued
fraction and exactly one even ordinary continued fraction. \qed
\end{theorem}

\subsection{Definition of lattice trigonometric functions.}

In this subsection we define the functions lattice sine, tangent,
and cosine on the set of ordinary lattice angles and formulate
their basic properties. We describe a geometric interpretation of
lattice trigonometric functions in terms of ordinary continued
fractions. Then we give the definitions of ordinary lattice
angles that are adjacent, transpose, and opposite interior to the
given angles. We use the notions of adjacent and transpose
ordinary lattice angles to define ordinary lattice right angles.

Let $A$, $B$, and $C$ be three lattice points that do not lie in
the same straight line. We denote the ordinary lattice angle with
the vertex at $B$ and the rays $BA$ and $BC$ by $\angle ABC$.

One can chose any other lattice point $B'$ in the open lattice
ray $AB$ (but not in $AC$) and any lattice point $C'$ in the open
lattice ray $AC$. For us the ordinary lattice angle $\angle AOB$
{\it coincides} with the ordinary lattice angle $\angle A'O'B'$.
Further we denote this by $\angle AOB = \angle A'O'B'$.

\begin{definition}
Two ordinary lattice angles $\angle AOB$ and $\angle A'O'B'$ are
said to be {\it lattice-congruent} if there exist a
lattice-affine transformation which takes the point $O$ to $O'$
and the rays $OA$ and $OB$ to the rays $O'A'$ and $O'B'$
respectively. We denote this as follows: $\angle AOB \cong \angle
A'O'B'$.
\end{definition}

Here we note that the relation $\angle AOB\cong \angle BOA$ holds
only for special ordinary lattice angles. (See below in
Subsubsection~\ref{sec_transpose}.)

\subsubsection{Definition of lattice sine,
tangent, and cosine for an ordinary lattice angle}

\begin{figure}[h]
$$\epsfbox{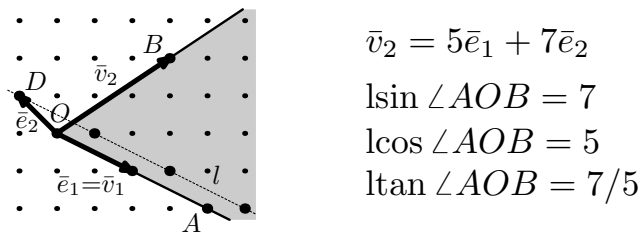}$$
\caption{An ordinary lattice angle $\angle AOB$ and its lattice
trigonometric functions.}\label{example7_5.10}
\end{figure}

Consider an arbitrary ordinary lattice angle $\angle AOB$. Let us
associate a special basis to this angle. Denote by $\bar v_1$ and
by $\bar v_2$ the lattice vectors generating the rays of the
angle:
$$
\bar v_1=\frac{\bar{OA}}{\il (OA)}, \quad \hbox{and} \quad \bar
v_2=\frac{\bar{OB}}{\il (OB)}.
$$
The set of lattice points on unit lattice distance from the
lattice straight line $OA$ coincides with the set of all lattice
points of two lattice straight lines parallel to $OA$. Since the
vectors $\bar v_1$ and $\bar v_2$ are linearly independent, the
ray $OB$ intersects exactly one of the above two lattice straight
lines. Denote this straight line by $l$. The intersection point
of the ray $OB$ with the straight line $l$ divides $l$ onto two
parts. Choose one of the parts which lies in the complement to
the convex hull of the union of the rays $OA$ and $OB$, and denote
by $D$ the lattice point closest to the intersection of the ray
$OB$ with the straight line $l$ (see Figure~\ref{example7_5.10}).

Now we choose the vectors $\bar e_1=\bar v_1$ and $\bar
e_2=\bar{OD}$. These two vectors are linearly independent and
generate the lattice. The basis $(\bar e_1,\bar e_2)$ is said to
be {\it associated} to the angle $\angle AOB$.

Since $(\bar e_1,\bar e_2)$ is a basis, the vector $\bar v_2$ has
a unique representation of the form:
$$
\bar v_2 = x_1 \bar e_1 +x_2 \bar e_2,
$$
where $x_1$ and $x_2$ are some integers.

\begin{definition}
In the above notation, the coordinates $x_2$ and $x_1$ are said to
be the {\it lattice sine} and the {\it lattice cosine} of the
ordinary lattice angle $\angle AOB$ respectively. The ratio of
the lattice sine and the lattice cosine ($x_2/x_1$) is said to be
the lattice tangent of $\angle AOB$.
\end{definition}

On Figure~\ref{example7_5.10} we show an example of lattice angle
with the lattice sine equals $7$ and the lattice cosine equals
$5$.

Let us briefly enumerate some elementary properties of lattice
trigonometric functions.

\begin{proposition}\label{basic_properties}
{\bf a).} The lattice sine and cosine of any ordinary lattice
angle are relatively-prime positive integers.
\\
{\bf b).} The values of lattice trigonometric functions for
lattice-congruent ordinary lattice angles coincide.
\\
{\bf c).} The lattice sine of an ordinary lattice angle coincide
with the index of the sublattice generated by all lattice vectors
of two angle rays in the lattice.
\\
{\bf d).} For any ordinary lattice angle $\alpha$ the following
inequalities hold:
$$
\isin\alpha\ge\icos\alpha, \quad \hbox{and} \quad
{\itan\alpha}\ge 1.
$$
The equalities hold iff the lattice vectors of the angle rays
generate the whole lattice.
\\
{\bf e). (Description of lattice angles.)} Two ordinary lattice
angles $\alpha$ and $\beta$ are lattice-congruent iff $\itan
\alpha =\itan \beta$. \qed
\end{proposition}

\subsubsection{Lattice arctangent.}
Let us fix the origin $O$ and a lattice basis $\bar e_1$ and
$\bar e_2$.
\begin{definition}
Consider an arbitrary rational $p\ge 1$. Let $p=m/n$, where $m$
and $n$ are positive integers. Suppose $A=O+\bar e_1$, and
$B=O+n\bar e_1 + m\bar e_2$. The ordinary angle $\angle AOB$ is
said to be the {\it arctangent of $p$ in the fixed basis} and
denoted by $\iarctan (p)$.
\end{definition}

The invariance of lattice tangents immediately implies the
following properties.

\begin{proposition}\label{itan_iarctan}
{\bf a).} For any rational $s\ge 1$, we have: $\itan(\iarctan
s)=s$.
\\
{\bf b).} For any ordinary lattice angle $\alpha$ the following
holds: $\iarctan(\itan \alpha) \cong \alpha$. \qed
\end{proposition}

\subsubsection{Lattice tangents, length-sine sequences, sails, and continued fractions.}
Let us start with the notion of sails for the ordinary lattice
angles. This notion is taken from theory of multidimensional
continued fractions in the sense of Klein (see, for example, the
works of F.~Klein~\cite{Kle1}, and
V.~Arnold~\cite{Arn2}).\label{KLEIN}

Consider an ordinary lattice angle $\angle AOB$. Let also the
vectors $\bar{OA}$ and $\bar{OB}$ be linearly independent, and of
unit lattice length.

Denote the closed convex solid cone for the ordinary lattice
angle $\angle AOB$ by $C(AOB)$. The boundary of the convex hull
of all lattice points of the cone $C(AOB)$ except the origin is
homeomorphic to the straight line. This boundary contains the
points $A$ and $B$. The closed part of this boundary contained
between the points $A$ and $B$ is called the {\it sail} for the
cone $C(AOB)$.

A lattice point of the sail is said to be a {\it vertex} of the
sail if there is no lattice segment of the sail containing this
point in the interior. The sail of the cone $C(AOB)$ is a broken
line with a finite number of vertices and without self
intersections. Let us orient the sail in the direction from $A$
to $B$, and denote the vertices of the sail by $V_i$ (for $0\le i
\le n$) according to the orientation of the sail (such that
$V_0=A$, and $V_n=B$).

\begin{definition}
Let the vectors $\bar{OA}$ and $\bar{OB}$ of the ordinary lattice
angle $\angle AOB$ be linearly independent, and of unit lattice
length. Let $V_i$, where $0\le i \le n$, be the vertices of the
corresponding sail. The sequence of lattice lengths and sines
$$
\begin{aligned}
(\il(V_0V_1),\isin\angle V_0V_1V_2,\il(V_1V_2),\isin\angle
V_1V_2V_3,
\ldots\\
\ldots,\il(V_{n-2}V_{n-1}),\isin\angle
V_{n-2}V_{n-1}V_n,\il(V_{n-1}V_n))
\end{aligned}
$$
is called the {\it lattice length-sine} sequence for the ordinary
lattice angle $\angle AOB$.
\end{definition}

\begin{remark}
The elements of the lattice length-sine sequence for any ordinary
lattice angle are positive integers. The lattice length-sine
sequences of lattice-congruent ordinary lattice angles coincide.
\end{remark}

\begin{theorem}\label{itan}
Let $ (\il(V_0V_1),\isin\angle V_0V_1V_2, \ldots, \isin\angle
V_{n-2}V_{n-1}V_n,\il(V_{n-1}V_n)) $ be the lattice length-sine
sequence for the ordinary lattice angle $\angle AOB$. Then the
{\it lattice tangent} of the ordinary lattice angle $\angle AOB$
equals to the value of the following ordinary continued fraction
$$
[\il(V_0V_1),\isin\angle V_0V_1V_2, \ldots, \isin\angle
V_{n-2}V_{n-1}V_n,\il(V_{n-1}V_n)].
$$
\end{theorem}

\begin{figure}[h]
$$\epsfbox{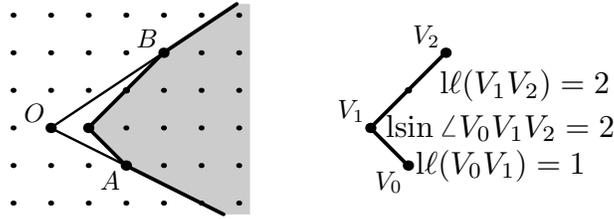}$$
\caption{ ${\itan\angle
AOB}=\frac{7}{5}=1+\frac{1}{2+1/2}$.}\label{example7_5}
\end{figure}

On Figure~\ref{example7_5} we show an example of an ordinary
lattice angle with tangent equivalent to $7/5$.

Further in Theorem~\ref{brokline} we formulate and prove a
general statement for generalized sails and signed length-sine
sequences. In the proof of Theorem~\ref{brokline} we refer only
on the preceding statements and definitions of
Subsection~\ref{expang_0}, that are independent of the statements
and theorems of all previous sections. For these reasons we skip
now the proof of Theorem~\ref{itan} (see also
Remark~\ref{proof_of_itan_iarctan}).

\subsubsection{Adjacent, transpose, and opposite interior ordinary lattice angles.}
\label{sec_transpose} First, we give the definition of the
ordinary lattice angles transpose and adjacent to the given one.

\begin{definition}
An ordinary lattice angle $\angle BOA$ is said to be {\it
transpose} to the ordinary lattice angle $\angle AOB$. We denote
it by $(\angle AOB)^{t}$.
\\
An ordinary lattice angle $\angle BOA'$ is said to be {\it
adjacent} to an ordinary lattice angle $\angle AOB$ if the points
$A$, $O$, and  $A'$ are contained in the same straight line, and
the point $O$ lies between $A$ and $A'$. We denote the ordinary
angle $\angle BOA'$ by $\pi {-}\angle AOB$.
\\
The ordinary lattice angle is said to be {\it right} if it is
self-dual and lattice-congruent to the adjacent ordinary angle.
\end{definition}

It immediately follows from the definition, that for any ordinary
lattice angle $\alpha$ the angles  $(\alpha^t)^t$ and $\pi
{-}(\pi{-}\alpha)$ are lattice-congruent to $\alpha$.

Further we will use the following notion. Suppose that some
integers $a$, $b$ and $c$, where $c\ge 1$, satisfy the following:
$ab\equiv 1 (\MOD c)$. Then we denote
$$
a\equiv \bigl(b (\MOD c)\bigr)^{-1}.
$$

For the trigonometric functions of transpose and adjacent ordinary
lattice angles the following relations are known.

\begin{theorem}\label{transpose_adjacent}
Consider an ordinary lattice angle $\alpha$. If $\alpha\cong
\iarctan (1)$, then
$$\alpha ^t\cong\pi{-}\alpha \cong \iarctan (1).$$
Suppose now, that $\alpha \not\cong \iarctan (1)$, then
$$
\begin{array}{c}
\isin(\alpha^t)=\isin\alpha,\quad
\icos(\alpha^t)\equiv \bigl(\icos\alpha(\MOD \isin\alpha)\bigr)^{-1};\\
\isin(\pi{-}\alpha)=\isin\alpha,\quad
\icos(\pi{-}\alpha)\equiv \bigl(-\icos\alpha(\MOD \isin\alpha)\bigr)^{-1}.\\
\end{array}
$$
Note also that $\pi{-}\alpha\cong \iarctan^t\left(\frac{\itan
\alpha}{\itan (\alpha)-1}\right)$. \qed
\end{theorem}

Theorem~\ref{transpose_adjacent} (after applying
Theorem~\ref{itan}) immediately reduces to the theorem of
P.~Popescu-Pampu. We refer the readers to his work~\cite{P-P} for
the proofs.

\subsubsection{Right ordinary lattice angles.}

It turns out that in lattice geometry there exist exactly two
lattice non-equivalent right ordinary lattice angles.

\begin{corollary}\label{right_}
Any ordinary lattice right angle is lattice-congruent to exactly
one of the following two angles: $\iarctan (1)$, or $\iarctan
(2)$. \qed
\end{corollary}

\begin{definition}
Consider two lattice parallel distinct straight lines $AB$ and
$CD$, where $A$, $B$, $C$, and $D$ are lattice points. Let the
points $A$ and $D$ be in different open half-planes with respect
to the straight line $BC$. Then the ordinary lattice angle
$\angle ABC$ is called {\it opposite interior} to the ordinary
lattice angle $\angle DCB$.
\end{definition}

Further we use the following proposition on opposite interior
ordinary lattice angles.

\begin{proposition}\label{opposite}
Two opposite interior to each other ordinary lattice angles are
lattice-congruent. \qed
\end{proposition}

The proof is left for the reader as an exercise.

\subsection{Basic lattice trigonometry of lattice angles in lattice triangles.}\label{srit}

In this subsection we introduce the sine formula for angles and
edges of lattice triangles. Further we show how to find the
lattice tangents of all angles and the lattice lengths of all
edges of any lattice triangle, if the lattice lengths of two
edges and the lattice tangent of the angle between them are given.

Let $A$, $B$, and $C$ be three distinct lattice points being not
contained in the same straight line. We denote the lattice
triangle with the vertices $A$, $B$, and $C$ by $\triangle ABC$.

\begin{definition}
The lattice triangles $\triangle ABC$ and $\triangle A'B'C'$ are
said to be {\it lattice-congruent} if there exist a
lattice-affine transformation which takes the point $A$ to $A'$,
$B$ to $B'$, and $C$ to $C'$ respectively. We denote: $\triangle
ABC {\cong} \triangle A'B'C'$.
\end{definition}

\subsubsection{The sine formula.}
Let us introduce a lattice analog of the Euclidean sine formula
for sines of angles and lengths of edges of triangles.

\begin{proposition}\label{sine}{\bf (The sine formula for lattice triangles.)}
The following holds for any lattice triangle $\triangle ABC$.
$$
\frac{\il (AB)}{\isin\angle BCA}= \frac{\il (BC)}{\isin\angle
CAB}= \frac{\il (CA)}{\isin\angle ABC}= \frac{\il (AB)\il (BC)\il
(CA)}{\is(\triangle ABC)}.
$$
\end{proposition}

\begin{proof}
The statement of Proposition~\ref{sine} follows directly from the
definition of lattice sine.
\end{proof}

\subsubsection{On relation between lattice tangents of ordinary lattice angles and lattice
lengths of lattice triangles.}

Suppose that we know the lattice lengths of the edges $AB$, $AC$
and the lattice tangent of $\angle BAC$ in the triangle
$\triangle ABC$. Now we show how to restore the lattice length
and the lattice tangents for the the remaining edge and ordinary
angles of the triangle.

For the simplicity we fix some lattice basis and use the system
of coordinates $OXY$ corresponding to this basis (denoted
$(*,*)$).

\begin{theorem}\label{ang-length}
Consider some triangle $\triangle ABC$. Let
$$
\il(AB)=c, \qquad \il(AC)=b, \quad \hbox{and} \quad \angle CAB
\cong \alpha.
$$
Then the ordinary angles $\angle BCA$ and $\angle ABC$ are
defined in the following way.
$$
\begin{array}{l}
\angle BCA \cong \left\{
\begin{array}{lcl}
\iarctan\left( \pi-\frac{c\isin \alpha}{c\icos\alpha -b} \right) &\hbox{if}&c\icos\alpha>b\\
\iarctan(1) &\hbox{if}&c\icos\alpha=b\\
\iarctan^t\left( \frac{c\isin \alpha}{b-c\icos\alpha} \right) &\hbox{if}&c\icos\alpha<b\\
\end{array}
\right.
,\\
\angle ABC \cong \left\{
\begin{array}{lcl}
\iarctan^t\left( \pi-\frac{b\isin (\alpha^t)}{b\icos(\alpha^t)-c} \right) &\hbox{if}&b\icos(\alpha^t)>c\\
\iarctan(1) &\hbox{if}&b\icos(\alpha^t)=c\\
\iarctan\left( \frac{b\isin (\alpha^t)}{c-b\icos(\alpha^t)} \right) &\hbox{if}&b\icos(\alpha^t)<c\\
\end{array}
\right.
.\\
\end{array}
$$
For the lattice length of the edge $CB$ we have
$$
\frac{\il(CB)}{\isin\alpha}=\frac{b}{\isin\angle
ABC}=\frac{c}{\isin\angle BCA}.
$$
\end{theorem}

\begin{proof}
Let $\alpha\cong\iarctan(p/q)$, where $\gcd(p,q)=1$. Then
$\triangle CAB \cong \triangle DOE$ where $D=(b,0)$, $O=(0,0)$,
and $E=(qc,pc)$. Let us now find the ordinary lattice angle
$\angle DEO$. Denote by $Q$ the point $(qc,0)$. If $qc{-}b= 0$,
then $\angle BCA=\angle DEO=\iarctan 1$. If $qc{-}b\ne 0$, then
we have
$$
\angle QDE \cong \iarctan\left(\frac{cp}{cq-b}\right)
 \cong\iarctan\left(\frac{c\isin \alpha}{c\icos\alpha-b}\right).
$$
The expression for $\angle BCA$ follows directly from the above
expression for $\angle QDE$, since $\angle BCA\cong\angle QDE$.
(See Figure~\ref{cases_tri}: here $\il(OD)=b$,
$\il(OQ)=c\icos\alpha$, and therefore $\il(DQ)=|c\icos\alpha-b|$.)
\begin{figure}[h]
$$\epsfbox{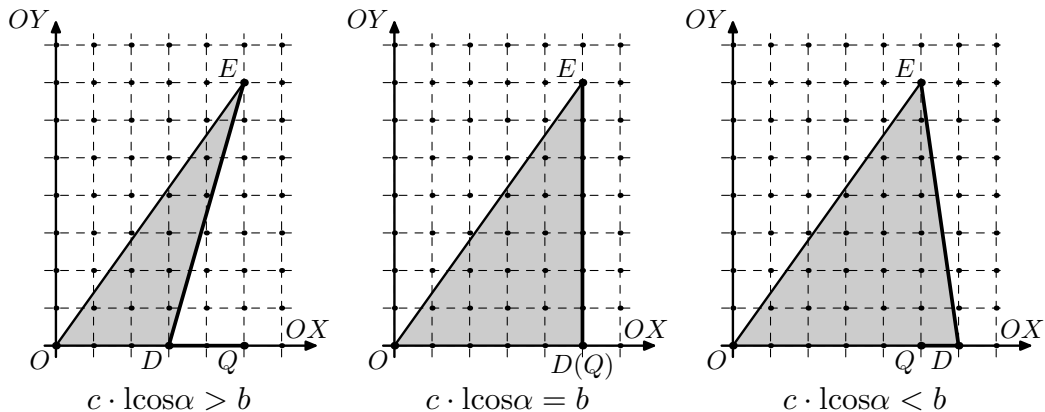}$$
\caption{Three possible configuration of points $O$, $D$, and
$Q$.}\label{cases_tri}
\end{figure}

To obtain the expression for $\angle ABC$ we consider the
triangle $\triangle BAC$. Calculate  $\angle CBA$ and then
transpose all ordinary angles in the expression. Since
$$
\begin{array}{l}
\is (ABC)=\il(AB)\il(AC)\isin\angle CAB=\il(BA)\il(BC)\isin\angle BCA=\\
\il(CB)\il(CA)\isin\angle ABC,
\end{array}
$$
we have the last statement of the theorem.
\end{proof}

\section{Theorem on sum of lattice tangents for the
ordinary lattice angles of lattice triangles. Proof of its second
statement.}

In this section we introduce the sum formula for the lattice
tangents of ordinary lattice angles of lattice triangles. The sum
formula is a lattice version of the following Euclidean statement:
three angles are the angles of some triangle iff their sum equals
$\pi$.

Throughout this section we fix some lattice basis and use the
system of coordinates $OXY$ corresponding to this basis.

\subsection{Finite continued fractions with not necessary positive elements.}\label{q1q2q3}

We start this section with the commonly-used definition of finite
continued fractions with not necessary positive elements. Let us
expand the set of rationals $\q$ with the operations $+$ and
$1/*$ on it with the element $\infty$. We pose
$q\pm\infty=\infty$, $1/0=\infty$, $1/\infty=0$ (we do not define
$\infty\pm\infty$ here). Denote this expansion by $\bar{\q}$.

For any finite sequence of integers $(a_0, a_1, \ldots ,a_n)$ we
associate an element $q$ of $\bar{\q}$:
$$
\begin{array}{ccc}
q&=&a_0+\frac{\displaystyle 1}{\displaystyle a_1+\frac{1}{
\begin{array}{cc}
\ddots&\vdots\\
& a_{n-1}+\frac{1}{a_n}\\
\end{array}
}}
\end{array},
$$
and denote it by $]a_0, a_1, \ldots ,a_n[$.

Let $q_i$ be some rationals, $i=1,\ldots,k$. Suppose that the odd
continued fraction for $q_i$ is $[a_{i,0}, a_{i,1}, \ldots
,a_{i,2n_i}]$ for $i=1,\ldots,k$. We denote by $]q_1, q_2, \ldots
,q_n[$ the following number
$$
\begin{aligned}
]a_{1,0}, a_{1;1}, \ldots ,a_{1,2n_1}, a_{2,0}, a_{2,1}, \ldots
,a_{2,2n_2},\ldots
a_{k-1,0}, a_{k-1,1}, \ldots ,a_{k-1,2n_{k-1}},\\
a_{k,0}, a_{k,1}, \ldots ,a_{k,2n_k}[
\end{aligned}
$$

\subsection{Formulation of the theorem and proof of its second statement.}

In Euclidean geometry the sum of Euclidean angles of the triangle
equals $\pi$. For any $3$-tuple of angles with the sum equals
$\pi$ there exist a triangle with these angles. Two Euclidean
triangles with the same angles are homothetic. Let us show one
generalization of these statements to the case of lattice
geometry.

We start this subsection with a definition of convex polygons
$n$-tuple to other convex polygons.

Let $n$ be an arbitrary positive integer, and $A=(x,y)$ be an
arbitrary lattice point. Denote by $nA$ the point $(nx,ny)$.
\begin{definition}
Consider any convex polygon or broken line with vertices $A_0,
\ldots, A_k$. The polygon or broken line $nA_0 \ldots nA_k$ is
called {\it n-multiple} (or {\it multiple}) to the given polygon
or broken line.
\end{definition}

\begin{theorem}\label{sum}{\bf On sum of lattice tangents of
ordinary lattice angles of lattice triangles.}\\
{\bf a).} Let $(\alpha_1,\alpha_2,\alpha_3)$ be an ordered
$3$-tuple of ordinary lattice angles. There exist a triangle with
three consecutive ordinary angles lattice-congruent to
$\alpha_1$, $\alpha_2$, and $\alpha_3$ iff there exist $i\in
\{1,2,3\}$ such that the angles $\alpha=\alpha_i$,
$\beta=\alpha_{i+1 (\MOD 3)}$, and $\gamma=\alpha_{i+2 (\MOD 3)}$
satisfy the following conditions:

{\it i$)$} the rational $]\itan\alpha ,-1,\itan\beta [$ is either
negative or greater than $\itan\alpha$;

{\it ii$)$} $]\itan\alpha ,-1,\itan\beta ,-1,\itan\gamma [=0$.
\\
{\bf b).} Let the consecutive ordinary lattice angles of some
triangle be $\alpha$, $\beta$, and $\gamma$. Then this triangle
is multiple to the triangle with vertices  $A_0=(0,0)$,
$B_0=(\lambda_2 \icos \alpha, \lambda_2 \isin \alpha)$, and
$C_0=(\lambda_1,0)$, where
$$
\lambda_1=\frac{\lcm(\isin\alpha, \isin\beta,
\isin\gamma)}{\gcd(\isin\alpha, \isin\gamma)}, \quad \hbox{and}
\quad \lambda_2=\frac{\lcm(\isin\alpha, \isin\beta,
\isin\gamma)}{\gcd(\isin\alpha, \isin\beta)}.
$$
\end{theorem}

\begin{remark}
Note that the statement of the Theorem~\ref{sum}a holds only for
odd continued fractions for the tangents of the correspondent
angles. We illustrate this with the following example. Consider a
lattice triangle with the lattice area equals $7$ and all angles
lattice-congruent to $\iarctan 7/3$ (see on Figure~\ref{triang}).
If we take the odd continued fractions $7/3=[2,2,1]$ for all
lattice angles of the triangle, then we have
$$
]2,2,1,-1,2,2,1,-1,2,2,1[=0.
$$
If we take the even continued fractions $7/3=[2,3]$ for all
angles of the triangle, then we have
$$
]2,3,-1,2,3,-1,2,3[=\frac{35}{13} \ne 0.
$$
\end{remark}

We prove the first statement of this theorem further in
Subsections~\ref{section_sum} and~\ref{section_sum2} after making
definitions of expanded lattice angles and their sums, and
introducing their properties.

{\it Proof of the second statement of Theorem~\ref{sum}.}
Consider the triangle $\triangle ABC$ with the ordinary lattice
angles $\alpha$, $\beta$, and $\gamma$ (with vertices at $A$,
$B$, and $C$ respectively). Suppose that for any $k>1$ and any
lattice triangle $\triangle KLM$ the triangle $\triangle ABC$ is
not lattice-congruent to the $k$-multiple of $\triangle KLM$.
That is equivalent to the following
$$
\gcd\bigl(\il(AB),\il(BC),\il(CA)\bigr)=1.
$$

Suppose that $S$ is the lattice area of $\triangle ABC$. Then by
the sine formula the following holds
$$
\left\{
\begin{array}{ccc}
\il(AB)\il(AC)&=&S/\isin\alpha\\
\il(BC)\il(BA)&=&S/\isin\beta\\
\il(CA)\il(CB)&=&S/\isin\gamma\\
\end{array}
\right. .
$$
Since $\gcd(\il(AB),\il(BC),\il(CA))=1$, we have
$\il(AB)=\lambda_1$ and $\il(AC)=\lambda_2$.

Therefore, the lattice triangle $\triangle ABC$ is
lattice-congruent to the lattice triangle $\triangle A_0B_0C_0$
of the theorem. \qed

\section{Expansion of ordinary lattice angles. The notions of sums
for lattice angles.}

In this section we introduce the notion of expanded lattice
angles and their normal forms and give the definition of sums of
expanded lattice angles. From the definition of  sums of
expanded  lattice angles we obtain a definition of  sums of
ordinary lattice angles. For the definition of expanded lattice
angles we expand the notion of sails in the sense of Klein: we
define and study oriented broken lines on the unit distance from
lattice points.

Throughout this section we work in the oriented two-dimensional
real vector space with a fixed lattice. We again fix some lattice
(positively oriented) basis and use the system of coordinates
$OXY$ corresponding to this basis.

The lattice-affine transformation is said to be {\it proper} if
it is orientation-preserving. We say that two sets are {\it
proper lattice-congruent} to each other if there exist a proper
lattice-affine transformation of $\r^2$ taking the first set to
the second.

\subsection{On a particular generalization of sails in the sense of Klein.}\label{expang_0}

In this subsection we introduce the definition of an oriented
broken lines on the unit lattice distance from a lattice point.
This notion is a direct generalization of the notion of s sail in
the sense of Klein (see page~\pageref{KLEIN} for the definition
of a sail). We expand the definition of lattice length-sine
sequences and continued fractions to the case of these broken
lines. We show that expanded lattice length-sine sequence for
oriented broken lines uniquely identifies the proper
lattice-congruence class of the corresponding broken line.
Further, we study the geometrical interpretation of the
corresponding continued fraction.

\subsubsection{Definition of a lattice signed length-sine sequence.}

For the definition of expanded lattice angles we expand the
definition of lattice length-sine sequence to the case of certain
broken lines.

For the 3-tuples of lattice points $A$, $B$, and $C$ we define
the function $\sgn$ as follows:
$$
\sgn(ABC)= \left\{
\begin{array}{ll}
+1,&\hbox{if the couple of vectors $\bar{BA}$ and $\bar{BC}$ defines the positive }\\
&\hbox{orientation.}\\
0,&\hbox{if the points $A$, $B$, and $C$ are contained in the same straight line.}\\
-1,&\hbox{if the couple of vectors $\bar{BA}$ and $\bar{BC}$ defines the negative}\\
&\hbox{orientation.}\\
\end{array}
\right.
$$

We also denote by $\sign:\r \to \{-1,0,1\}$ the sign function
over reals.

The segment $AB$ is said to be {\it on the unit distance} from
the point $C$ if the lattice vectors of the segment $AB$, and the
vector $\bar{AC}$ generate the lattice.

A union of (ordered) lattice segments $A_0A_1, A_1A_2, \ldots,
A_{n-1}A_n$ ($n>0$) is said to be a {\it lattice oriented broken
line} and denoted by $A_0A_1A_2 \ldots A_n$ if any two
consecutive segments are not contained in the same straight line.
We also say that the lattice oriented broken line
$A_nA_{n-1}A_{n-2} \ldots A_0$ is {\it inverse} to the lattice
oriented broken line $A_0A_1A_2 \ldots A_n$.

\begin{definition}
Consider a lattice oriented broken line and a lattice point in
the complement to this line. The broken line is said to be {\it
on the unit distance} from the point if all edges of the broken
line are on the unit lattice distance from the point.
\end{definition}

Let us now associate to any lattice oriented broken line on the
unit distance from some point the following sequence of non-zero
elements.

\begin{definition}
Let $A_0A_1\ldots A_n$ be a lattice oriented broken line on the
unit distance from some lattice point $V$. The sequence $(a_0,
\ldots, a_{2n-2})$ defined as follows:
$$
\begin{array}{l}
a_0=\sgn(A_0VA_1)\il(A_0A_1),\\
a_1=\sgn(A_0VA_1)\sgn(A_1VA_2)\sgn(A_0A_1A_2)\isin\angle A_0A_1A_2,\\
a_2=\sgn(A_1VA_2)\il(A_1A_2),\\
\ldots\\
a_{2n-3}=\sgn(A_{n-2}VA_{n-1})\sgn(A_{n-1}VA_{n})\sgn(A_{n-2}A_{n-1}A_n)\isin\angle A_{n-2}A_{n-1}A_{n},\\
a_{2n-2}=\sgn(A_{n-1}VA_n)\il(A_{n-1}A_n),
\end{array}
$$
is called an {\it lattice signed length-sine} sequence for the
lattice oriented broken line on the unit distance
from the lattice point $V$.\\
The element of $\bar{\q}$
$$
]a_0,a_1,\ldots,a_{2n-2}[
$$
is called the {\it continued fraction for the broken line}
$A_0A_1\ldots A_n$.
\end{definition}

\begin{figure}[h]
$$\epsfbox{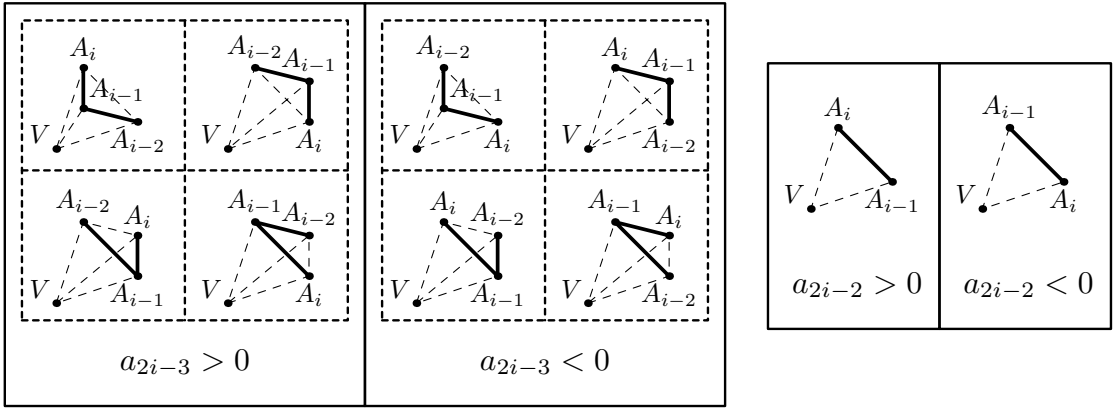}$$
\caption{All possible (non-degenerate) proper affine
decompositions for angles and segments of a signed length-sine
sequence.}\label{signs}
\end{figure}

We show how to identify signs of elements of the lattice signed
length-sine sequence for a lattice oriented broken line on the
unit distance from the lattice point $V$ on Figure~\ref{signs}.

\begin{figure}[h]
$$\epsfbox{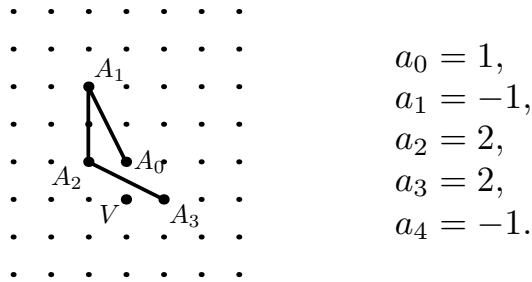}$$
\caption{A lattice oriented broken line on the unit distance from
the point $V$ and the corresponding lattice signed length-sine
sequence.}\label{m1_2first}
\end{figure}

On Figure~\ref{m1_2first} we show an example of lattice oriented
broken line on the unit distance from the lattice point $V$ and
the corresponding lattice signed length-sine sequence.

\begin{proposition}\label{ccc}
A lattice signed length-sine sequence for the given lattice
oriented broken line and the lattice point is invariant under the
group action of the proper lattice-affine transformations.
\end{proposition}

\begin{proof}
The statement of the proposition holds, since the functions
$\sgn$, $\il$, and $\isin$ are invariant under the group action
of the proper lattice-affine transformations.
\end{proof}

\subsubsection{Proper lattice-congruence of lattice oriented broken lines
on the unit distance from lattice points.}

Let us formulate necessary and sufficient conditions for two
lattice oriented broken lines on the unit distance from the same
lattice point to be proper lattice-congruent.

\begin{theorem}\label{twobrok}
The lattice signed length-sine sequences of two lattice oriented
broken lines on the unit distance from lattice points $V_1$ and
$V_2$ coincide iff there exist a proper lattice-affine
transformation taking the point $V_1$ to $V_2$ and one lattice
oriented broken line to the other.
\end{theorem}

\begin{proof}
The lattice signed length-sine sequence for any lattice oriented
broken line on the unit distance is uniquely defined, and by
Proposition~\ref{ccc} is invariant under the group action of
lattice-affine orientation preserving transformations. Therefore,
the lattice signed length-sine sequences for two proper
lattice-congruent lattice oriented broken lines coincide.

\vspace{2mm}

Suppose now that two lattice oriented broken lines $A_0 \ldots
A_n$, and $B_0 \ldots B_n$ on the unit distance from lattice
points $V_1$ and $V_2$ respectively have the same lattice signed
length-sine sequence  $(a_0, a_1, \ldots , a_{2n-3}, a_{2n-2})$.
Let us prove that these broken lines are proper lattice-congruent.
Without loose of generality we consider the point $V_1$ at the
origin $O$.

Let $\xi$ be the proper lattice-affine transformation taking the
point $V_2$ to the point $V_1=O$, $B_0$ to $A_0$, and the lattice
straight line containing $B_0B_1$ to the lattice straight line
containing $A_0A_1$. Let us prove inductively that $\xi
(B_i)=A_i$.

\vspace{2mm}

{\it Base of induction.} Since $a_0=b_0$, we have
$$
\sgn(A_0OA_1)\il(A_0A_1)=\sgn(\xi(B_0)O\xi(B_1))\il(\xi(B_0)\xi(B_1)).
$$
Thus, the lattice segments $A_0A_1$ and $A_0\xi(B_1)$ are of the
same lattice length and of the same direction. Therefore,
$\xi(B_1)=A_1$.

\vspace{2mm}

{\it Step of induction.} Suppose, that $\xi (B_i)=A_i$ holds for
any nonnegative integer $i\le k$, where $k\ge 1$. Let us prove,
that $\xi(B_{k+1})=A_{k+1}$. Denote by $C_{k+1}$ the lattice
point $\xi(B_{k+1})$. Let $A_{k}=(q_k,p_k)$. Denote by $A_{k}'$
the closest lattice point of the segment $A_{k-1}A_k$ to the
vertex $A_k$. Suppose that $A_k'=(q_k',p_k')$. We know also
$$
\begin{array}{lll}
a_{2k-1}&=&\sgn(A_{k-1}OA_{k})\sgn(A_kOC_{k+1})\sgn(A_{k-1}A_kC_{k+1})\isin\angle A_{k-1}A_kC_{k+1},\\
a_{2k}&=&\sgn(A_kOC_{k+1})\il(A_kC_{k+1}).
\end{array}
$$

Let the coordinates of $C_{k+1}$ be $(x,y)$. Since
$\il(A_kC_{k+1})=|a_{2k}|$ and the segment $A_kC_{k+1}$ is on the
unit lattice distance to the origin $O$, we have $\is(\triangle
OA_kC_{k+1})=|a_{2k}|$. Since the segment $OA_k$ is of the unit
lattice length, the coordinates of $C_{k+1}$ satisfy the
following equation:
$$
|-p_kx+q_ky|=|a_{2k}|.
$$
Since $\sgn(A_kOC_{k+1})\il(A_kC_{k+1})=\sign(a_{2k})$, we have
$$
-p_kx+q_ky=a_{2k}.
$$
Since $\isin\angle A_{k}'A_kC_{k+1}=\isin\angle
A_{k-1}A_kC_{k+1}=|a_{2k-1}|$, and the lattice lengths of
$A_kC_{k+1}$, and $A_{k}'A_k$ are $|a_{2k}|$ and $1$ respectively,
we have $\is(\triangle A_k'A_kC_{k+1})=|a_{2k-1}a_{2k}|$.
Therefore, the coordinates of $C_{k+1}$ satisfy the following
equation:
$$
|-(p_k-p_k')(x-q_k)+(q_k-q_k')(y-p_k)|=|a_{2k-1}a_{2k}|.
$$
Since
$$
\left\{
\begin{array}{rcl}
\sgn(A_{k-1}OA_{k})\sgn(A_kOC_{k+1})\sgn(A_{k-1}A_kC_{k+1})&=&\sign(a_{2k-1})\\
\sgn(A_kOC_{k+1})&=&\sign(a_{2k})
\end{array}
\right. ,
$$
we have
$$
(p_k-p_k')(x-q_k)-(q_k-q_k')(y-p_k)=\sgn(A_{k-1}OA_{k})a_{2k-1}a_{2k}.
$$

We obtain the following:
$$
\left\{
\begin{array}{rcl}
-p_kx+q_ky&=&a_{2k}\\
(p_k-p_k')(x-q_k)-(q_k-q_k')(y-p_k)&=&\sgn(A_{k-1}OA_{k})a_{2k-1}a_{2k}\\
\end{array}
\right. .
$$
Since
$$
\left|\det \left(
\begin{array}{cc}
-p_k& q_k\\
p_k'-p_k& q_k-q_k'\\
\end{array}
\right) \right| =1,
$$
there exist and unique an integer solution for the system of
equations for $x$ and $y$. Hence, the points $A_{k+1}$ and
$C_{k+1}$ have the same coordinates. Therefore,
$\xi(B_{k+1})=A_{k+1}$. We have proven the step of induction.

The proof of Theorem~\ref{twobrok} is completed by induction.
\end{proof}

\subsubsection{Values of continued fractions for
lattice oriented broken lines on the unit distance from the
origin.}

Now we show the relation between lattice oriented broken lines on
the unit distance from the origin and the corresponding continued
fractions for them.

\begin{theorem}\label{brokline}
Let $A_0 A_1 \ldots A_n$ be a lattice oriented broken line on the
unit distance from the origin $O$. Let also $A_0=(1,0)$,
$A_1=(1,a_0)$, $A_n=(p,q)$, where $\gcd(p,q)=1$, and $(a_0, a_1,
\ldots ,a_{2n-2})$ be the corresponding lattice signed
length-sine sequence. Then the following holds:
$$
\frac{q}{p}=]a_0, a_1, \ldots ,a_{2n-2}[.
$$
\end{theorem}

\begin{proof}
To prove this theorem we use an induction on the number of edges
of the broken lines.

\vspace{2mm}

{\it Base of induction.} Suppose that a lattice oriented broken
line on the unit distance from the origin has a unique edge, and
the corresponding sequence is $(a_0)$. Then $A_1=(1,a_0)$ by the
assumptions of the theorem. Therefore, we have
$\frac{a_0}{1}{=}]a_0[$.

\vspace{2mm}

{\it Step of induction.} Suppose that the statement of the
theorem is correct for any lattice oriented broken line on the
unit distance from the origin with $k$ edges. Let us prove the
theorem for the arbitrary lattice oriented broken line on the
unit distance from the origin with $k{+}1$ edges (and satisfying
the conditions of the theorem).

Let $A_0  \ldots A_{k+1}$ be a lattice oriented broken line on
the unit distance from the origin with the lattice signed
length-sine sequence  $(a_0, a_1, \ldots , a_{2k-1}, a_{2k})$.
Let also
$$
A_0=(1,0), \qquad A_1=(1,a_0),\quad \hbox{and} \quad
A_{k+1}=(p,q).
$$
Consider the lattice oriented broken line $B_1 \ldots B_{k+1}$ on
the unit distance from the origin with a lattice signed
length-sine sequence  $(a_2, a_3, \ldots , a_{2k-2}, a_{2k})$.
Let also
$$
B_1=(1,0), \qquad B_2=(1,a_2), \quad \hbox{and} \quad
B_{k+1}=(p',q').
$$
By the induction assumption we have
$$
\frac{q'}{p'}=]a_2, a_3, \ldots ,a_{2k}[.
$$

We extend the lattice oriented broken line $B_1 \ldots B_{k+1}$
to the lattice oriented broken line $B_0 B_1  \ldots B_{k+1}$ on
the unit distance from the origin, where $B_0=(1{+}a_0a_1,-a_0)$.
Let the lattice signed length-sine sequence for this broken line
be $(b_0, b_1, \ldots , b_{2k-1}, b_{2k})$. Note that
$$
\begin{array}{lll}
b_0&=&\sgn(B_0OB_1)\il(B_0B_1)=\sign(a_0)|a_0|=a_0,\\
b_1&=&\sgn(B_0OB_1)\sgn(B_1OB_2)\sgn(B_0B_1B_2)\isin\angle B_0B_1B_2=\\
&&\sign{a_0}\sign{b_2}\sign(a_0a_1b_2)|a_1|=a_1,\\
b_l&=&a_l, \quad \hbox{for $l=2,\ldots, 2k$}.\\
\end{array}
$$

Consider a proper lattice-linear transformation $\xi$ that takes
the point $B_0$ to the point $(1,0)$, and $B_1$ to $(1,a_0)$.
These two conditions uniquely defines $\xi$:
$$
\xi= \left(
\begin{array}{cc}
1&a_1\\
a_0&1+a_0a_1\\
\end{array}
\right).
$$
Since $B_{k+1}=(p',q')$, we have
$\xi(B_{k+1})=(p'{+}a_1q',q'a_0{+}p'{+}p'a_0a_1)$.
$$
\frac{q'a_1+p'+p'a_0a_1}{p'+a_1q'}=a_0+\frac{1}{a_1+q'/p'}=]a_0,a_1,a_2,a_3,
\ldots ,a_{2n}[.
$$

Since, by Theorem~\ref{twobrok} the lattice oriented broken lines
$B_0 B_1 \ldots B_{k+1}$ and $A_0 A_1\ldots A_{k+1}$ are
lattice-linear equivalent, $B_0=A_0$, and $B_1=A_1$, these broken
lines coincide. Therefore, for the coordinates $(p,q)$ the
following hold
$$
\frac{q}{p}=\frac{q'a_0+p'+p'a_0a_1}{p'+a_1q'}=]a_0,a_1,a_2,a_3,
\ldots ,a_{2k}[.
$$

On Figure~\ref{m1_2second} we illustrate the step of induction
with an example of lattice oriented broken line on the unit
distance from the origin with the lattice signed length-sine
sequence: $(1,-1,2,2,-1)$. We start (the left picture) from the
broken line $B_1B_2B_3$ with the lattice signed length-sine
sequence: $(2,2,-1)$. Then, (the picture in the middle) we expand
the broken line $B_1B_2B_3$ to the broken line $B_0B_1B_2B_3$
with the lattice signed length-sine sequence: $(1,-1,2,2,-1)$.
Finally (the right picture) we apply a corresponding proper
lattice-linear transformation $\xi$ to achieve the resulting
broken line $A_0A_1A_2A_3$.
\begin{figure}[h]
$$\epsfbox{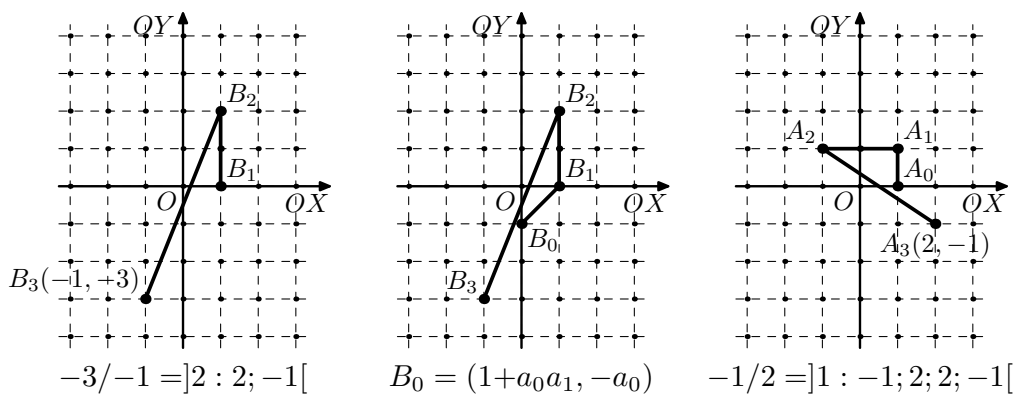}$$
\caption{The case of lattice oriented broken line on the unit
distance from the origin with lattice signed length-sine
sequence: $(1,-1,2,2,-1)$.}\label{m1_2second}
\end{figure}

We have proven the step of induction.

The proof of Theorem~\ref{brokline} is completed.

\end{proof}

\begin{remark}\label{proof_of_itan_iarctan}
Theorem~\ref{brokline} immediately implies the statement of
Theorem~\ref{itan}. One should put the sail of an angle as an
oriented-broken line $A_0A_1\ldots A_n$.
\end{remark}

\subsection{Expanded lattice angles. Sums for ordinary and expanded lattice angles.}\label{expang}

In this subsection we introduce the notion of expanded lattice
angles and their normal forms. We use normal forms of expanded
lattice angles to give the definition of sums of expanded lattice
angles. From the definition of  sums of expanded  lattice angles
we obtain a definition of  sums of ordinary  lattice angles.

\subsubsection{Equivalence classes of lattice oriented broken lines and
the corresponding expanded lattice angles.}

Let us give a definition of expanded  lattice angles.

\begin{definition}
Two lattice oriented broken lines $l_1$ and $l_2$ on the unit
distance from $V$ are said to be {\it equivalent} if they have in
common the first and the last vertices and the closed broken line
generated by $l_1$ and the inverse of $l_2$ is homotopy
equivalent to the point in $\r^2\setminus \{V\}$.\\
An equivalence class of lattice oriented broken lines on the unit
distance from $V$ containing the broken line $A_0 A_1 \ldots A_n$
is called the {\it expanded  lattice angle for the equivalence
class of $A_0A_1\ldots A_n$} at the {\it vertex} $V$ (or, for
short, {\it expanded  lattice angle}) and denoted by $\angle
(V,A_0A_1 \ldots A_n)$.
\end{definition}

We study the expanded lattice angles up to proper
lattice-congruence (and not up to lattice-congruence).

\begin{definition}
Two expanded  lattice angles $\Phi_1$ and $\Phi_2$ are said to be
{\it proper lattice-congruent} iff there exist a proper
lattice-affine transformation sending the class of lattice
oriented broken lines corresponding to $\Phi_1$ to the class of
lattice oriented broken lines corresponding to $\Phi_2$. We
denote this by $\Phi_1\pcong \Phi_2$.
\end{definition}

\subsubsection{Revolution numbers for expanded  lattice angles.}

Let us describe one invariant of expanded  lattice angles under
the group action of the proper lattice-affine transformations.

Let $r=\{V{+}\lambda\bar{v}|\lambda\ge 0\}$ be the oriented ray
for an arbitrary vector $\bar v$ with the vertex at $V$, and $AB$
be an oriented (from $A$ to $B$) segment not contained in the ray
$r$. Suppose also, that the vertex $V$ of the ray $r$ is not
contained in the segment $AB$. We denote by $\# (r,V,AB)$ the
following number:
$$
\# (r,V,AB)= \left\{
\begin{array}{ll}
0,& \hbox{if the segment $AB$ does not intersect the ray $r$}\\
\frac{\displaystyle 1}{\displaystyle 2}
\sgn\Bigl(A(A{+}\bar{v}) B\Bigr),& \hbox{if the segment $AB$ intersects the ray $r$ at $A$ or}\\
&\quad \hbox{ at $B$}\\
\sgn\Bigl(A(A{+}\bar{v}) B\Bigr),& \hbox{if the segment $AB$ intersects the ray $r$ at the}\\
&\quad \hbox{ interior point of $AB$}\\
\end{array}
\right. ,
$$
and call it the {\it intersection number} of the ray $r$ and the
segment $AB$.

\begin{definition}
Let $A_0A_1\ldots A_n$ be some lattice oriented broken line, and
let $r$ be an oriented ray $\{V{+}\lambda\bar{v}|\lambda\ge 0\}$.
Suppose that the ray $r$ does not contain the edges of the broken
line, and the broken line does not contain the point $V$. We call
the number
$$
\sum\limits_{i=1}^{n}\# (r,V,A_{i-1}A_i)
$$
the {\it intersection number} of the ray $r$ and the lattice
oriented broken line $A_0A_1\ldots A_n$, and denote it by $\#
(r,V,A_0A_1\ldots A_n)$.
\end{definition}

\begin{definition}
Consider an arbitrary expanded  lattice angle $\angle (V,
A_0A_1\ldots A_n)$. Denote the rays
$\{V+\lambda\bar{VA_0}|\lambda\ge 0\}$ and
$\{V-\lambda\bar{VA_0}|\lambda\ge 0\}$ by $r_+$ and $r_-$
respectively. The number
$$
\frac{1}{2}\bigl(\#(r_+,V,A_0A_1\ldots A_n)+\#(r_-,V,A_0A_1\ldots
A_n)\bigr)
$$
is called the {\it lattice revolution number} for the expanded
lattice angle $\angle (V,A_0A_1\ldots A_n)$, and denoted by
$\#(\angle (V,A_0A_1\ldots A_n)).$
\end{definition}

Let us prove that the above definition is correct.

\begin{proposition}
The revolution number of any expanded  lattice angle is
well-defined.
\end{proposition}

\begin{proof}
Consider an arbitrary expanded  lattice angle $\angle
(V,A_0A_1\ldots A_n)$. Let
$$
r_+=\{V+\lambda\bar{VA_0}|\lambda\ge 0\} \quad \hbox{ and} \quad
r_-=\{V-\lambda\bar{VA_0}|\lambda\ge 0\}.
$$

Since the lattice oriented broken line $A_0A_1\ldots A_n$ is on
the unit lattice distance from the point $V$, any segment of this
broken line is on the unit lattice distance from $V$. Thus, the
broken line does not contain $V$, and the rays $r_+$ and $r_-$ do
not contain edges of the curve.

Suppose that
$$
\angle (V,A_0A_1\ldots A_n)=\angle (V',A'_0A_1'\ldots A_m').
$$ \
This implies that $V=V'$, $A_0=A_0'$, $A_n=A_m'$, and the broken
line
$$
A_0A_1\ldots A_nA_{m-1}'\ldots A_1'A_0'
$$
is homotopy equivalent to the point in $\r^2\setminus \{V\}$.
Thus,
$$
\begin{array}{l}
\#(\angle (V,A_0A_1\ldots A_n))-\#(\angle (V,'A'_0A_1'\ldots A_m'))=\\
\frac{1}{2}\bigl(\#(r_+,V,A_0A_1\ldots A_nA_{m-1}'\ldots A_1'
A_0')+
\#(r_-,V,A_0A_1\ldots A_nA_{m-1}'\ldots A_1' A_0')\bigr)=\\
0{+}0=0.
\end{array}
$$
Hence,
$$
\#(\angle (V,A_0A_1\ldots A_n))=\#(\angle (V',A'_0A_1\ldots
A_m')).
$$
Therefore, the revolution number of any expanded  lattice angle
is well-defined.
\end{proof}

\begin{proposition}
The revolution number of expanded lattice angles is invariant
under the group action of the proper lattice-affine
transformations. \qed
\end{proposition}

\subsubsection{Zero ordinary angles.}
For the next theorem we will need to define define zero ordinary
angles and their trigonometric functions. Let $A$, $B$, and $C$
be three lattice points of the same lattice straight line.
Suppose that $B$ is distinct to $A$ and $C$ and the rays $BA$ and
$BC$ coincide. We say that the ordinary lattice angle with the
vertex at $B$ and the rays $BA$ and $BC$ is {\it zero}. Suppose
$\angle ABC$ is zero, put by definition
$$
\isin(\angle ABC)=0, \quad \icos(\angle ABC)=1, \quad
\itan(\angle ABC)=0.
$$
Denote by $\iarctan (0)$ the angle $\angle AOA$ where $A=(1,0)$,
and $O$ is the origin.

\subsubsection{On normal forms of expanded  lattice angles.}

Let us formulate and prove a theorem on normal forms of expanded
lattice angles. We use the following notation.

By the sequence
$$
\bigl((a_0,\ldots,a_n)\times
\hbox{$k$-times},b_0,\ldots,b_m\bigr),
$$
where $k\ge 0$, we denote the following sequence:
$$
(\underbrace{a_0,\ldots,a_n,a_0,\ldots,a_n, \quad \ldots\quad
,a_0,\ldots,a_n}_{\mbox{$k$-times}} ,b_0,\ldots,b_m).
$$

\begin{definition}\label{normal_forms}
$\bf I).$ Suppose $O$ be the origin, $A_0$ be the point $(1,0)$.
We say that the expanded lattice angle $\angle(O,A_0)$ is {\it of
the type} $\bf I$ and denote it by $0\pi+\iarctan(0)$ (or $0$,
for short). The empty sequence is said to be {\it characteristic}
for the angle $0\pi+\iarctan(0)$.
\vspace{1mm}\\
Consider a lattice oriented broken line $A_0A_1\ldots A_s$ on the
unit distance from the origin $O$. Let also $A_0$ be the point
$(1,0)$, and the point $A_1$ be on the straight line $x=1$. If
the signed length-sine sequence of the expanded ordinary angle
$\Phi_0=\angle(O,A_0A_1\ldots A_s)$
coincides with the following sequence (we call it {\it characteristic sequence} for the corresponding angle):\\
{\bf II$_k$)} $\bigl((1,-2,1,-2)\times
\hbox{$(k-1)$-times},1,-2,1\bigr)$, where $k\ge 1$, then we
denote the angle $\Phi_0$ by $k\pi{+}\iarctan(0)$ (or $k\pi$, for
short)
and say that $\Phi_0$ is {\it of the type} {\bf II$_k$};\\
{\bf III$_k$)} $\bigl((-1,2,-1,2)\times
\hbox{$(k-1)$-times},-1,2,-1\bigr)$, where $k\ge 1$, then we
denote the angle $\Phi_0$ by $-k\pi {+}\iarctan(0)$ (or $-k\pi$,
for short)
and say that $\Phi_0$ is {\it of the type} {\bf III$_k$};\\
{\bf IV$_k$)} $\bigl((1,-2,1,-2)\times \hbox{$k$-times},
a_0,\ldots, a_{2n}\bigr)$, where $k\ge 0$, $n\ge 0$, $a_i>0$, for
$i=0,\ldots, 2n$, then we denote the angle $\Phi_0$ by $k\pi
{+}\iarctan([a_0,a_1,\ldots, a_{2n}])$
and say that $\Phi_0$ is {\it of the type} {\bf IV$_k$};\\
{\bf V$_k$)} $\bigl((-1,2,-1,2)\times \hbox{$k$-times},
a_0,\ldots, a_{2n}\bigr)$, where $k> 0$, $n\ge 0$, $a_i>0$, for
$i=0,\ldots, 2n$, then we denote the angle $\Phi_0$ by $-k\pi {+}
\iarctan([a_0,a_1,\ldots, a_{2n}])$ and say that $\Phi_0$ is {\it
of the type} {\bf V$_k$}.
\end{definition}

\begin{theorem}\label{nf}
For any expanded  lattice angle $\Phi$ there exist and unique a
type among the types {\bf{I}}-{\bf{V}} and a unique expanded
lattice angle $\Phi_0$ of that type such that $\Phi_0$ is proper
lattice
congruent to $\Phi$.\\
{\rm The expanded  lattice angle $\Phi_0$ is said to be {\it the
normal form} for the expanded  lattice angle $\Phi$.}
\end{theorem}

For the proof of Theorem~\ref{nf} we need the following lemma.

\begin{lemma}\label{lema_nf}
Let $m$, $k\ge 1$, and $a_i>0$ for $i=0,\ldots, 2n$ be some integers.\\
{\bf a).} Suppose the lattice signed length-sine sequences for
the expanded  lattice angles $\Phi_1$ and $\Phi_2$ are
respectively
$$
\begin{array}{l}
\bigl((1,-2,1,-2)\times \hbox{$(k{-}1)$-times},1,-2,1,-2, a_0,\ldots, a_{2n}\bigr) \quad \hbox{and}\\
\bigl((1,-2,1,-2)\times \hbox{$(k{-}1)$-times},1,-2,1,m, a_0,\ldots, a_{2n}\bigr),\\
\end{array}
$$
then $\Phi_1$ is proper lattice-congruent to $\Phi_2$.\\
{\bf b).} Suppose the lattice signed length-sine sequences for
the expanded  lattice angles $\Phi_1$ and $\Phi_2$ are
respectively
$$
\begin{array}{l}
\bigl((-1,2,-1,2)\times \hbox{$(k{-}1)$-times},-1,2,-1,m, a_0,\ldots, a_{2n}\bigr) \quad \hbox{and}\\
\bigl((-1,2,-1,2)\times \hbox{$(k{-}1)$-times},-1,2,-1,2, a_0,\ldots, a_{2n}\bigr),\\
\end{array}
$$
then $\Phi_1$ is proper lattice-congruent to $\Phi_2$.
\end{lemma}

\begin{proof}
We prove the first statement of the lemma. Suppose that $m$ is
integer, $k$ is positive integer, and $a_i$ for $i=0,\ldots, 2n$
are positive integers.

Let us construct the angle $\Psi_1$ with vertex at the origin for
the lattice oriented broken line $A_0 \dots A_{2k+n+1}$,
corresponding to the lattice signed length-sine sequence
$$
\bigl((1,-2,1,-2)\times \hbox{$(k{-}1)$-times},1,-2,1,-2,
a_0,\ldots, a_{2n}\bigr),
$$
such that $A_0=(1,0)$, $A_1=(1,1)$. Note that
$$
\left\{
\begin{array}{llll}
A_{2l}&=&((-1)^l,0),&\hbox{for $l{<}k{-}1$}\\
A_{2l+1}&=&((-1)^l,(-1)^l) &\hbox{for $l{<}k{-}1$}\\
A_{2k}&=&((-1)^k,0)\\
A_{2k+1}&=&((-1)^k,(-1)^k)a_0\\
\end{array}
\right. .
$$

Let us construct the angle $\Psi_2$ with vertex at the origin for
the lattice oriented broken line $B_0 \dots B_{2k+n+1}$,
corresponding to the lattice signed length-sine sequence
$$
\bigl((1,-2,1,-2)\times \hbox{$(k{-}1)$-times},1,-2,1,m,
a_0,\ldots, a_{2n}\bigr).
$$
such that $B_0=(1,0)$, $B_1=(-m-1,1)$. Note also that
$$
\left\{
\begin{array}{llll}
B_{2l}&=&((-1)^l,0),&\hbox{for $l{<}k{-}1$}\\
B_{2l+1}&=&((-1)^l(-m-1),(-1)^l) &\hbox{for $l{<}k{-}1$}\\
B_{2k}&=&((-1)^k,0)\\
B_{2k+1}&=&((-1)^k,(-1)^k)a_0\\
\end{array}
\right. .
$$

From the above we know, that the points $A_{2k}$ and $A_{2k+1}$
coincide with the points $B_{2k}$ and $B_{2k+1}$ respectively.
Since the remaining parts of both lattice signed length-sine
sequences (i.~e. $(a_0,\ldots, a_{2n})$) coincide, the point
$A_l$ coincide with the point $B_l$ for $l>2k$.

Since the lattice oriented broken lines $A_0  \ldots  A_{2k}$ and
$B_0 \ldots B_{2k}$ are of the same equivalence class, and  the
point $A_l$ coincide with the point $B_l$ for $l>2k$, we obtain
$$
\Psi_1=\angle(O,A_0 \ldots A_{2k+n+1})=\angle(O,B_0 \ldots
B_{2k+n+1})=\Psi_2.
$$

Therefore, by Theorem~\ref{twobrok} we have the following:
$$
\Phi_1\pcong\Psi_1=\Psi_2\pcong\Phi_2.
$$
This concludes the proof of Lemma~\ref{lema_nf}a.

Since the proof of Lemma~\ref{lema_nf}b almost completely repeats
the proof of Lemma~\ref{lema_nf}a, we omit the proof of
Lemma~\ref{lema_nf}b here.
\end{proof}

{\it Proof of Theorem~\ref{nf}.} First, we prove that any two
distinct expanded  lattice angles listed in
Definition~\ref{normal_forms} are not proper lattice-congruent.
Let us note that the revolution numbers of expanded  lattice
angles distinguish the types of the angles. The revolution number
for the expanded  lattice angle of the type {\bf I} is $0$. The
revolution number for the expanded  lattice angle of the type
{\bf II$_k$} is $1/2(k{+}1)$ where $k\ge 0$. The revolution
number for the expanded  lattice angle of the type {\bf III$_k$}
is $-1/2(k{+}1)$ where $k\ge 0$. The revolution number for the
expanded  lattice angles of the type {\bf IV$_k$} is $1/4{+}1/2k$
where $k\ge 0$. The revolution number for the expanded  lattice
angles of the type {\bf V$_k$} is $1/4{-}1/2k$ where $k > 0$.

So we have proven that two expanded  lattice angles of different
types are not proper lattice-congruent. For the types {\bf I},
{\bf II$_k$}, and {\bf III$_k$} the proof is completed, since any
such type consists of the unique expanded  lattice angle.

Let us prove that normal forms of the same type {\bf IV$_k$} (or
of the same type {\bf V$_k$}) are not proper lattice-congruent
for any integer $k\ge 0$ (or $k>0$). Consider an expanded
lattice angle $\Phi=k\pi {+} \iarctan([a_0,a_1,\ldots, a_{2n}])$.
Suppose that a lattice oriented broken line  $A_0A_1 \ldots A_m$
on the unit distance from $O$, where $m=2|k|{+}n{+}1$ defines the
angle $\Phi$. Suppose also that the signed lattice-sine sequence
for this broken line is characteristic.

Suppose, that $k$ is even, then the ordinary  lattice angle
$\angle A_0OA_m$ is proper lattice-congruent to the ordinary
lattice angle $\iarctan([a_0,a_1,\ldots, a_{2n}])$. This angle is
a proper lattice-affine invariant for the expanded  lattice angle
$\Phi$. This invariant distinguish the expanded  lattice angles
of type {\bf IV$_k$} (or {\bf V$_k$}) with even $k$.

Suppose, that $k$ is odd, then denote $B=O{+}\bar{A_0O}$. The
ordinary  lattice angle $\angle BVA_m$ is proper
lattice-congruent to the ordinary  lattice angle
 $\iarctan([a_0,a_1,\ldots, a_{2n}])$.
This angle is a proper lattice-affine invariant for the expanded
lattice angle $\Phi$. This invariant distinguish the expanded
lattice angles of type {\bf IV$_k$} (or {\bf V$_k$}) with odd $k$.

Therefore, the expanded  lattice angles listed in
Definition~\ref{normal_forms} are not proper lattice-congruent.

\vspace{2mm}

Now we prove that an arbitrary expanded  lattice angle is proper
lattice-congruent to one of the expanded  lattice angles of the
types {\bf{I}}-{\bf{V}}.

Consider an arbitrary expanded  lattice angle $\angle
(V,A_0A_1\ldots A_n)$ and denote it by $\Phi$. If $\#(\Phi)=k/2$
for some integer $k$, then $\Phi$ is proper lattice congruent to
an angle of one of the types {\bf{I}}-{\bf{III}}. Let
$\#(\Phi)=1/4$, then the expanded  lattice angle $\Phi$ is proper
lattice-congruent to the expanded  lattice angle defined by the
sail of the ordinary  lattice angle $\angle A_0VA_n $ of the type
{\bf IV$_0$}.

Suppose now, that $\#(\Phi)=1/4{+}k/2$ for some positive integer
$k$, then one of its lattice signed length-sine sequence is of
the following form:
$$
\bigl((1,-2,1,-2)\times \hbox{$(k-1)$-times},1,-2,1,m,
a_0,\ldots, a_{2n}\bigr),
$$
where $a_i>0$, for $i=0,\ldots, 2n$. By Lemma~\ref{lema_nf} the
expanded  lattice angle defined by this sequence is proper
lattice-congruent to an expanded  lattice angle of the type {\bf
IV$_k$} defined by the sequence
$$
\bigl((1,-2,1,-2)\times \hbox{$(k-1)$-times},1,-2,1,-2,
a_0,\ldots, a_{2n}\bigr).
$$

Finally, let $\#(\Phi)=1/4{-}k/2$ for some positive integer $k$,
then one of its lattice signed length-sine sequence is of the
following form:
$$
\bigl((-1,2,-1,2)\times \hbox{$(k-1)$-times},-1,2,-1,m,
a_0,\ldots, a_{2n}\bigr),
$$
where $a_i>0$, for $i=0,\ldots, 2n$. By Lemma~\ref{lema_nf} the
expanded  lattice angle defined by this sequence is proper
lattice-congruent to an expanded  lattice angle of the type {\bf
V$_k$} defined by the sequence
$$
\bigl((-1,2,-1,2)\times \hbox{$(k-1)$-times},-1,2,-1,2,
a_0,\ldots, a_{2n}\bigr).
$$

This completes the proof of Theorem~\ref{nf}. \qed

Let us finally give the definition of trigonometric functions for
the expanded lattice angles and describe some relations between
ordinary and expanded lattice angles.

\begin{definition}
Consider an arbitrary expanded  lattice angle $\Phi$ with the
normal form $k\pi {+} \varphi$ for some ordinary (possible zero)
lattice angle
$\varphi$ and for an integer $k$.\\
{\bf a).} The ordinary  lattice angle $\varphi$ is said to be
{\it associated} with the expanded  lattice angle $\Phi$.\\
{\bf b).} The numbers $\itan(\varphi)$, $\isin(\varphi)$, and
$\icos(\varphi)$ are called the {\it lattice tangent}, the {\it
lattice sine}, and the {\it lattice cosine} of the expanded
lattice angle $\Phi$.
\end{definition}

Since all sails for ordinary  lattice angles are lattice oriented
broken lines, the set of all ordinary angles is naturally
embedded into the set of expanded  lattice angles.

\begin{definition}
For any ordinary  lattice angle $\varphi$ the angle
$$
0\pi + \iarctan(\itan\varphi)
$$
is said to be {\it corresponding} to the angle $\varphi$ and
denoted by $\bar\varphi$.
\end{definition}

From Theorem~\ref{nf} it follows that for any ordinary  lattice
angle $\varphi$ there exists and unique an expanded  lattice
angle $\bar \varphi$ corresponding to $\varphi$. Therefore, two
ordinary lattice angles $\varphi_1$ and $\varphi_2$ are
lattice-congruent iff the corresponding  lattice angles $\bar
\varphi_1$ and $\bar \varphi_2$ are proper lattice-congruent.

\subsubsection{Opposite expanded  lattice angles. Sums of expanded  lattice angles.
Sums of ordinary  lattice angles.} Consider an expanded  lattice
angle $\Phi$ with the vertex $V$ for some equivalence class of a
given lattice oriented broken line. The expanded lattice angle
$\Psi$ with the vertex $V$ for the equivalence class of the
inverse lattice oriented broken line is called {\it opposite} to
the given one and denoted by $-\Phi$.

\begin{proposition}
For any expanded  lattice angle $\Phi\pcong k\pi  {+}\varphi$ we
have:
$$
-\Phi\pcong (-k-1)\pi +(\pi-\varphi).
$$
\qed
\end{proposition}

Let us introduce the definition of sums of ordinary and expanded
lattice angles.

\begin{definition}
Consider arbitrary expanded  lattice angles $\Phi_i$,
$i=1,\ldots,l$. Let the characteristic sequences for the normal
forms of $\Phi_i$ be $(a_{0,i},a_{1,i},\ldots, a_{2n_i,i})$ for
$i=1,\ldots,l$. Let $M=(m_1,\ldots,m_{l-1}$) be some
$(l{-}1)$-tuple of integers. The normal form of any expanded
lattice angle, corresponding to the following lattice signed
length-sine sequence
$$
\begin{aligned}
\bigl(a_{0,1},a_{1,1},\ldots, a_{2n_1,1},m_1,
a_{0,2},a_{1,2},\ldots, a_{2n_2,2},m_2, \ldots \\
\ldots ,m_{l-1},a_{0,l},a_{1,l},\ldots, a_{2n_l,l}\bigr),
\end{aligned}
$$
is called the {\it $M$-sum of expanded  lattice angles $\Phi_i$
$($$i=1,\ldots,l$$)$} and denoted by
$$
\sum\limits_{M,i=1}^{l}\Phi_i, \quad \hbox{or equivalently by}
\quad \Phi_1 +_{m_1} \Phi_2+_{m_2} \ldots+_{m_{l-1}}\Phi_l.
$$
\end{definition}

\begin{proposition}
The $M$-sum of expanded  lattice angles $\Phi_i$
$($$i=1,\ldots,l$$)$ is well-defined. \qed
\end{proposition}

Let us say a few words about properties of $M$-sums.

The $M$-sum of expanded  lattice angles is non-associative. For
example, let $\Phi_1\pcong  \iarctan 2$, $\Phi_2\pcong \iarctan
(3/2)$, and $\Phi_3\pcong \iarctan 5$. Then
$$
\begin{array}{rcl}
\Phi_1+_{-1}\Phi_2+_{-1}\Phi_3&=& \pi + \iarctan(4),\\
\Phi_1+_{-1}(\Phi_2+_{-1}\Phi_3)&=& 2\pi,\\
(\Phi_1+_{-1}\Phi_2)+_{-1}\Phi_3&=& \iarctan(1).\\
\end{array}
$$

The $M$-sum of expanded  lattice angles is non-commutative. For
example, let $\Phi_1\pcong \iarctan 1$, and $\Phi_2\pcong
\iarctan 5/2$. Then
$$
\begin{array}{rcl}
\Phi_1+_{1}\Phi_2&=& \iarctan(12/7),\\
\Phi_2+_{1}\Phi_1&=& \iarctan(13/5).\\
\end{array}
$$

\begin{remark}
The $M$-sum of expanded  lattice angles is naturally extended to
the sum of classes of proper lattice-congruences of expanded
lattice angles.
\end{remark}

We conclude this section with the definition of sums of ordinary
lattice angles.

\begin{definition}
Consider ordinary lattice angles $\alpha_i$, where $i=1,\ldots,l$.
Let $\bar\alpha_i$ be the corresponding expanded lattice angles
for $\alpha_i$, and $M=(m_1,\ldots,m_{l-1}$) be some
$(l{-}1)$-tuple of integers. The ordinary  lattice angle $\varphi$
associated with the expanded  lattice angle
$$
\Phi=\bar\alpha_1 +_{m_1}
\bar\alpha_2+_{m_2}\ldots+_{m_{l-1}}\bar\alpha_l.
$$
is called the {\it $M$-sum of ordinary  lattice angles $\alpha_i$
$($$i=1,\ldots,l$$)$} and denoted by
$$
\sum\limits_{M,i=1}^{l}\alpha_i, \quad \hbox{or equivalently by}
\quad \alpha_1 +_{m_1} \alpha_2+_{m_2}\ldots+_{m_{l-1}}\alpha_l.
$$
\end{definition}

\begin{remark}
Note that the sum of ordinary  lattice angles is naturally
extended to the classes of lattice-congruences of  lattice angles.
\end{remark}

\section{Relations between expanded lattice angles and ordinary lattice angles.
Proof of the first statement of Theorem~\ref{sum}.}

In this section we show how to calculate the ordinary angle
$\varphi$ of the normal form: we describe relations between
continued fractions for lattice oriented broken lines and the
lattice tangents for the corresponding expanded lattice angles.
Then we prove the first statement of the theorem on sums of
lattice tangents for ordinary lattice angles in lattice
triangles. Further we define lattice-acute-angled triangles. We
conclude this section with a necessary and sufficient condition
for an ordered $n$-tuple of angles to be the angles of some
convex lattice polygon.

Throughout this section we again fix some lattice basis and use
the system of coordinates $OXY$ corresponding to this basis.

\subsection{On relations between continued fractions for lattice oriented broken lines
and the lattice tangents of the corresponding expanded lattice
angles.}

For any real number $r$ we denote by $\lfloor r \rfloor$ the
maximal integer not greater than $r$.

\begin{theorem}\label{phi}
Consider an expanded lattice angle $\Phi=\angle (V,A_0A_1\ldots
A_n)$. Suppose, that the normal form for $\Phi$ is $k\pi {+}
\varphi$ for some integer $k$ and an ordinary lattice angle
$\varphi$. Let $(a_0, a_1, \ldots ,a_{2n-2})$ be the lattice
signed length-sine sequence for the lattice oriented broken line
$A_0A_1 \ldots A_n$. Suppose that
$$
]a_0, a_1, \ldots ,a_{2n-2}[=q/p.
$$
Then the following holds:
$$
\varphi\cong \left\{
\begin{array}{ll}
\iarctan (1),&\hbox{if $q/p{=}\infty$,}\\
\iarctan(q/p),&\hbox{if $q/p{\ge}1$,}\\
\iarctan\left(\frac{|q|}{|p|-\lfloor (|p|-1)/|q| \rfloor |q|}\right), &\hbox{if $0{<}q/p{<}1$,}\\
0,&\hbox{if $q/p=0$,}\\
\pi-\iarctan\left(\frac{|q|}{|p|-\lfloor (|p|-1)/|q| \rfloor |q|}\right), &\hbox{if $-1{<}q/p{<}0$,}\\
\pi-\iarctan(-q/p),&\hbox{if $q/p{\le}-1$.}\\
\end{array}
\right.
$$
\end{theorem}

\begin{proof}
Consider the following linear coordinates $(*,*)'$ on the plane
$\r^2$, associated with the lattice oriented broken line $A_0 A_1
\ldots A_n$ and the point $V$. Let the origin $O'$ be at the
vertex $V$, $(1,0)'=A_0$, and
$(1,1)'=A_0+\frac{1}{a_0}\sgn(A_0O'A_1)\bar{A_0A_1}$. The other
coordinates are uniquely defined by linearity. We denote this
system of coordinates by $O'X'Y'$.

The set of integer points for the coordinate system $O'X'Y'$
coincides with the set of lattice points of $\r^2$. The basis of
vectors $(1,0)'$ and $(0,1)'$ defines a positive orientation.

Suppose that the new coordinates of the point $A_n$ are
$(p',q')'$. Then by Theorem~\ref{brokline} we have $q'/p'=q/p$.
This directly implies the statement of the theorem for the cases
$q'{>}p'{>}0$, $q'/p'=0$, and $q'/p'=\infty$.

Suppose now that $p'{>}q'{>}0$. Consider the ordinary lattice
angle $\varphi= \angle A_0PA_n$. Let $B_0\ldots B_m$ be the sail
for it. The direct calculations show that the point
$$
D=B_0+\frac{\bar{B_0B_1}}{\il (B_0B_1)}
$$
coincides with the point $(1{+}\lfloor (p'-1)/q'\rfloor,1)$ in
the system of coordinates $O'X'Y'$.

Consider the proper lattice-linear (in the coordinates $O'X'Y'$)
transformation $\xi$ that takes the point $A_0=B_0$ to itself,
and the point $D$ to the point $(1,1)'$. These conditions
uniquely identify $\xi$.
$$
\xi= \left(
\begin{array}{cc}
1& -\lfloor (p'-1)/q'\rfloor\\
0& 1\\
\end{array}
\right)
$$

The transformation $\xi$ takes the point $A_n=B_m$ with the
coordinates $(p',q')$ to the point with the coordinates
$(p'-\lfloor (p'-1)/q'\rfloor q',q')'$. Since $q'/p'=q/p$, we
obtain the following
$$
\varphi=\iarctan\left(\frac{q'}{p'-\lfloor (p'-1)/q'\rfloor
q'}\right)= \iarctan\left(\frac{q}{p-\lfloor (p-1)/q\rfloor
q}\right).
$$

\vspace{2mm}

The proof for the case $q'>0$ and $p'<0$ repeats the described
cases after taking to the consideration the adjacent angles.

Finally, the case of $q'<0$ repeats all previous cases by the
central symmetry (centered at the point $O'$) reasons.

This completes the proof of Theorem~\ref{phi}.
\end{proof}

\begin{corollary}
The revolution number and the continued fraction for any lattice
oriented broken line on the unit distance from the vertex
uniquely define the proper lattice-congruence class of the
corresponding expanded lattice angle. \qed
\end{corollary}

\subsection{Proof of Theorem~\ref{sum}a: two preliminary lemmas.}\label{section_sum}
We say that the lattice point $P$ is {\it on the lattice
distance} $k$ from the lattice segment $AB$ if the lattice
vectors of the segment $AB$ and the vector $\bar{AP}$ generate a
sublattice of the lattice of index $k$.

\begin{figure}[h]
$$\epsfbox{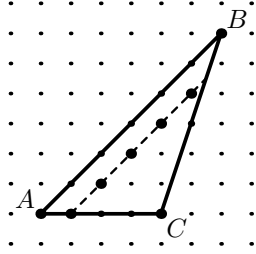}$$
\caption{For the given triangle $\triangle ABC$ we have
$\il_1(AB;C)=5$.}\label{sum3}
\end{figure}

\begin{definition}
Consider a lattice triangle $\triangle ABC$. Denote the number of
lattice points on the unit lattice distance from the segment $AB$
and contained in the $($closed$)$ triangle $\triangle ABC$ by
$\il_1(AB;C)$ (see on Figure~\ref{sum3}).
\end{definition}

Note that all lattice points on the lattice unit distance from the
segment $AB$ in the $($closed$)$ lattice triangle $\triangle ABC$
are contained in one straight line parallel to the straight line
$AB$. Besides, the integer $\il_1(AB;C)$ is positive for any
triangle $\triangle ABC$.

Now we prove the following lemma.

\begin{lemma}\label{formula}
For any lattice triangle $\triangle ABC$ the following holds
$$
\bar{\angle CAB}+_{\il(AB)-\il_1(AB;C)-1}\bar{\angle
ABC}+_{\il(BC)-\il_1(BC;A)-1}\bar{\angle BCA}=\pi .
$$
\end{lemma}

\begin{proof}
Consider an arbitrary lattice triangle $\triangle ABC$. Suppose
that the couple of vectors $\bar{BA}$ and $\bar{BC}$ defines the
positive orientation of the plane (otherwise we apply to the
triangle $\triangle ABC$ some lattice-affine transformation
changing the orientation and come to the same position). Denote
(see Figure~\ref{sum2} below):
$$
\begin{array}{l}
D=A+\bar{BC}, \quad \hbox{and} \quad E=A+\bar{AC}.\\
\end{array}
$$

Since $CADB$ is a parallelogram, the triangle $\triangle BAD$ is
proper lattice-congruent to the triangle $\triangle ABC$. Thus,
the angle $\angle BAD $ is proper lattice-congruent to the angle
$\angle ABC$, and $\il_1(BA;D)=\il_1(AB;C)$. Since $EABD$ is a
parallelogram, the triangle $\triangle AED$ is proper
lattice-congruent to the triangle $\triangle BAD$, and hence is
proper lattice-congruent to the triangle  $\triangle ABC$. Thus,
$\angle DAE$ is proper lattice-congruent to $\angle BCA$, and
$\il_1(DA;E)=\il_1(BC;A)$.

Let $A_0\ldots A_n$ be the sail of $\angle CAB$ with the
corresponding lattice length-sine sequence $(a_0,\ldots,
a_{2n-2})$. Let $B_0 B_1 \ldots B_m$ be the sail of $\angle BAD$
(where $B_0=A_n$) with the corresponding lattice length-sine
sequence $(b_0,\ldots, b_{2m-2})$. And let $C_0C_1 \ldots C_l$ be
the sail of $\angle DAE$ (where $C_0=B_m$) with the corresponding
lattice length-sine sequence $(c_0,\ldots, c_{2l-2})$.

Consider now the lattice oriented broken line
$$
A_0 \ldots A_nB_0 B_1 \ldots B_mC_0 C_1 \ldots C_l.
$$
The lattice oriented length-sine sequence for this broken line is
$$
(a_0,\ldots, a_{2n-2},t,b_0,\ldots, b_{2m-2},u,c_0,\ldots,
c_{2l-2}).
$$
By definition of the sum of expanded lattice angles this sequence
defines the expanded lattice angle
$$
\bar{\angle CAB}+_{t}\bar{\angle BAD}+_{u}\bar{\angle DAE}.
$$
By Theorem~\ref{phi},
$$
\bar{\angle CAB}+_t\bar{\angle BAD}+_u\bar{\angle DAE}= \pi.
$$

\vspace{2mm}

Let us find an integer $t$. Denote by $A_n'$ the closest lattice
point to the point $A_n$ and distinct to $A_n$ in the segment
$A_{n-1}A_n$. Consider the set of lattice points on the unit
distance from the segment $AB$ and lying in the half-plane with
the boundary straight line $AB$ and containing the point $D$.
This set coincides with the following set (See Figure~\ref{sum1}):
$$
\left\{A_{n,k}=A_n+\bar{A_n'A_n}+k\bar{AA_n}|k\in\z \right\}.
$$
\begin{figure}[h]
$$\epsfbox{sum1.1}$$
\caption{Lattice points $A_{n,t}$.}\label{sum1}
\end{figure}

Since $A_{n,-2}=A+\bar{A_n'A}$, the points $A_{n,k}$ for $k\le
-2$ are in the closed half-plane bounded by the straight line
$AC$ and not containing the point $B$.

Since $A_{n,-1}=A+\bar{A_n'A_n}$, the points $A_{n,k}$ for
$k\ge-1$ are in the open half-plane bounded by the straight line
$AC$ and containing the point $B$.

The intersection of the parallelogram $AEDB$ and the open
half-plane bounded by the straight line $AC$ and containing the
point $B$ contains exactly $\il (AB)$ points of the described set:
only the points $A_{n,k}$ with $-1\le k \le \il(AB){-}2$.

Since the triangle $\triangle BAD$ is proper lattice-congruent to
$\triangle ABC$, the number of points $A_{n,k}$ in the closed
triangle $\triangle BAD$ is $\il_1(AB;C)$: the points $A_{n,k}$
for
$$
\il(AB)-\il_1(AB;C)-1\le k \le \il(AB)-2.
$$
Denote the integer $\il(AB){-}\il_1(AB;C){-}1$ by $k_0$.

The point $A_{n,k_0}$ is contained in the segment $B_0B_1$ of the
sail for the ordinary lattice angle $\angle BAD$ (see
Figure~\ref{sum2}). Since the angles $\angle BAD$ and $\angle
ABC$ are proper lattice-congruent, we have
$$
\begin{aligned}
t=\sgn(A_{n-1}AA_n)\sgn(A_nAB_1)\sgn(A_{n-1}A_nB_1)\isin\angle A_{n-1}A_nB_1=\\
1\cdot 1\cdot \sgn(A_{n-1}A_nA_{n,k_0})\isin\angle A_{n-1}A_nA_{n,k_0}=\\
\sign (k_0)|k_0|=k_0=\il(AB)-\il_1(AB;C)-1.
\end{aligned}
$$
\begin{figure}[h]
$$\epsfbox{sum1.2}$$
\caption{The point $A_{n,k_0}$.}\label{sum2}
\end{figure}

Exactly by the same reasons,
$$
u=\il(DA)-\il_1(DA;E)-1=\il(BC)-\il_1(BC;A)-1.
$$
\vspace{2mm}

Therefore, $ \bar{\angle CAB}+_{\il(AB)-\il_1(AB;C)-1}\bar{\angle
ABC}+_{\il(BC)-\il_1(BC;A)-1}\bar{\angle BCA}=\pi.$
\end{proof}

\begin{lemma}\label{existence}
Let $\alpha$, $\beta$, and $\gamma$ be nonzero ordinary lattice
angles. Suppose that $\bar\alpha+_u\bar\beta+_v\bar\gamma=\pi$,
then there exist a triangle with three consecutive ordinary
angles lattice-congruent to $\alpha$, $\beta$, and~$\gamma$.
\end{lemma}

\begin{proof}
Denote by $O$ the point $(0,0)$, by $A$ the point $(1,0)$, and by
$D$ the point $(-1,0)$ in the fixed system of coordinates $OXY$.

Let us choose the points $B=(p_1,q_1)$ and $C=(p_2,q_2)$ with
integers $p_1$, $p_2$ and positive integers $q_1$, $q_2$ such that
$$
\angle AOB=\iarctan(\itan\alpha), \quad  \hbox{and}\quad \angle
AOC=\bar{\angle AOB}+_u\bar\beta.
$$
Thus the vectors $\bar{OB}$ and $\bar{OC}$ defines the positive
orientation, and $\angle BOC \cong \beta.$ Since
$$
\bar\alpha+_u\bar\beta+_v\bar\gamma=\pi \quad \hbox{and} \quad
\bar\alpha+_u\bar\beta \pcong \angle AOC,
$$
the ordinary angle $\angle COD$ is lattice-congruent to $\gamma$.

Denote by $B'$ the point $(p_1q_2,q_1q_2)$, and by $C'$ the point
$(p_2q_1,q_1q_2)$ and consider the triangle $B'OC'$. Since the
ordinary angle $\angle B'OC'$ coincides with the ordinary angle
$\angle BOC$, we obtain
$$
\angle B'OC'\cong\beta.
$$

Since the ordinary angle $\beta$ is nonzero, the points $B'$ and
$C'$ are distinct and the straight line $B'C'$ does not coincide
with the straight line $OA$. Since the second coordinate of the
both points $B'$ and $C'$ equal $q_1q_2$, the straight line
$B'C'$ is parallel to the straight line $OA$. Thus, by
Proposition~\ref{opposite} it follows that
$$
\angle C'B'O\cong \angle AOB'=\angle AOB\cong \alpha, \quad
\hbox{and} \quad \angle OC'B'\cong \angle C'OD=\angle COD\cong
\gamma.
$$

So, we have constructed the triangle $\triangle B'OC'$ with three
consecutive ordinary angles lattice-congruent to $\alpha$,
$\beta$, and $\gamma$.
\end{proof}

\subsection{Proof of Theorem~\ref{sum}a: conclusion of the proof.}\label{section_sum2}

Now we return to the proof of the first statement of the theorem
on sums of lattice tangents for ordinary lattice angles in
lattice triangles.

{\it Proof of Theorem~\ref{sum}a.} Let $\alpha$, $\beta$, and
$\gamma$ be nonzero ordinary lattice angles satisfying the
conditions {\it i$)$} and {\it ii$)$} of Theorem~\ref{sum}a.

The second condition
$$
]\itan(\alpha),-1,\itan(\beta),-1,\itan(\gamma)[=0
$$
implies that
$$
\bar{\alpha}+_{-1}\bar{\beta}+_{-1}\bar{\gamma}=k\pi.
$$
Since all three tangents are positive, we have $k=1$, or $k=2$.

Consider the first condition: $]\itan\alpha ,-1,\itan\beta [$ is
either negative or greater than $\itan\alpha$. It implies that
$\bar{\alpha}+_{-1}\bar{\beta}=0\pi+\varphi$, for some ordinary
lattice angle $\varphi$, and hence $k=1$.

Therefore, by Lemma~\ref{existence} there exist a triangle with
three consecutive ordinary lattice angles lattice-congruent to
$\alpha$, $\beta$, and $\gamma$.

\vspace{1mm}

Let us prove the converse. We prove that condition {\it ii$)$} of
Theorem~\ref{sum}a holds by reductio ad absurdum. Suppose, that
there exist a triangle $\triangle ABC$ with consecutive ordinary
angles $\alpha=\angle CAB$, $\beta=\angle ABC$, and
$\gamma=\angle BCA$, such that
$$
\left\{
\begin{array}{lcl}
]\itan(\alpha),-1,\itan(\beta),-1,\itan(\gamma)[&\ne &0\\
]\itan(\beta),-1,\itan(\gamma),-1,\itan(\alpha)[&\ne &0\\
]\itan(\gamma),-1,\itan(\alpha),-1,\itan(\beta)[&\ne &0\\
\end{array}
\right. .
$$

These inequalities and Lemma~\ref{formula} imply that at least two
of the integers
$$
\il(AB){-}\il_1(AB;C){-}1, \qquad \il(BC){-}\il_1(BC;A){-}1,
\quad \hbox{and} \quad \il(CA){-}\il_1(CA;B){-}1
$$
are nonnegative.

Without losses of generality we suppose that
$$
\left\{
\begin{array}{lcl}
\il(AB)-\il_1(AB;C)-1&\ge &0\\
\il(BC)-\il_1(BC;A)-1&\ge &0\\
\end{array}
\right. .
$$
Since all integers of the continued fraction
$$
r=]\itan(\alpha),\il(AB)-\il_1(AB;C)-1,\itan(\beta),\il(BC)-\il_1(BC;A)-1,\itan(\gamma)[
$$
are non-negative and the last one is positive, we obtain that
$r>0$ (or $r=\infty$). From the other hand, by
Lemma~\ref{formula} and by Theorem~\ref{phi} we have that
$r=0/{-1}=0$. We come to the contradiction.

Now we prove that condition {\it i$)$} of Theorem~\ref{sum}a
holds. Suppose that there exist a triangle $\triangle ABC$ with
consecutive ordinary angles $\alpha=\angle CAB$, $\beta=\angle
ABC$, and $\gamma=\angle BCA$, such that
$$
]\itan(\alpha),-1,\itan(\beta),-1,\itan(\gamma)[=0.
$$
Since $ \bar{\alpha}+_{-1}\bar{\beta}+_{-1}\bar{\gamma}=\pi, $ we
have $ \bar{\alpha}+_{-1}\bar{\beta}=0\pi+\varphi $ for some
ordinary lattice angle $\varphi$. Therefore, the first condition
of the theorem holds.

This concludes the proof of Theorem~\ref{sum}. \qed

\vspace{2mm}

Let us give here the following natural definition.
\begin{definition}
The triangle $\triangle ABC$ is said to be {\it
lattice-acute-angled} if the integers
$$
\il(AB){-}\il_1(AB;C){-}1, \qquad \il(BC){-}\il_1(BC;A){-}1,
\quad \hbox{and} \quad \il(CA){-}\il_1(CA;B){-}1
$$
are all equal to $-1$.\\
The triangle $\triangle ABC$ is said {\it lattice right-angled} if one of these integers equals $0$.\\
The triangle $\triangle ABC$ is said {\it lattice obtuse-angled}
if one of these integers is positive.
\end{definition}

Note that the property of the lattice triangle to be lattice
acute-angled, right-angled, or obtuse-angled cannot be determined
by one of the angles, unlike in Euclidean geometry. We will
illustrate this with the following example. (Actually, there is
no relation between right-angled triangles and right angles,
defined before in Subsubsection~\ref{sec_transpose}.)
\begin{example}
On Figure~\ref{obtuse} we show the lattice obtuse-angled triangle
(on the figure to the left) with the lattice tangents of its
ordinary angles equal to $3/2$, $8/3$, and $1$. Nevertheless,
each of these lattice angles is not supposed to be
``lattice-obtuse'', since three lattice triangles (on the figure
to the right) are all lattice acute-angled and contain ordinary
lattice angles with the lattice tangents $3/2$, $8/3$, and~$1$.
\begin{figure}[h]
$$\epsfbox{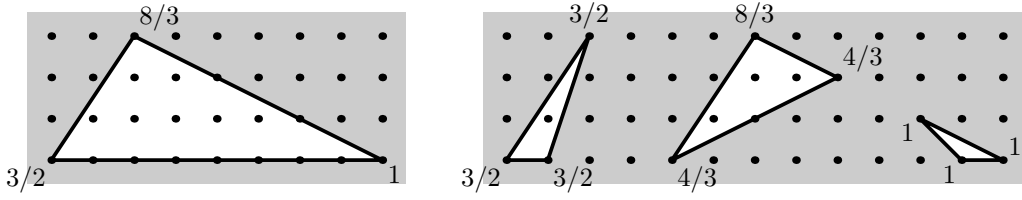}$$
\caption{The triangle to the left is lattice obtuse-angled, three
triangles to the right are all lattice
acute-angled.}\label{obtuse}
\end{figure}
\end{example}

\subsection{Theorem on sum of lattice tangents for ordinary lattice angles of convex polygons.}

A satisfactory description for lattice-congruence classes of
lattice convex polygons has not been yet found. It is only known
that the number of convex polygons with lattice area bounded from
above by $n$ growths exponentially in $n$, while $n$ tends to
infinity (see~\cite{Arn5} and~\cite{Bar2}). We conclude this
section with the following theorem on necessary and sufficient
condition for the lattice angles to be the angles of some convex
lattice polygon.

\begin{theorem}\label{sum_for_polygons}
Let $\alpha_1,\ldots, \alpha_n$ be an arbitrary ordered $n$-tuple
of ordinary non-zero lattice angles.
Then the following two conditions are equivalent:\\
-- there exist a convex $n$-vertex polygon with consecutive
ordinary lattice angles lattice-congruent
to the ordinary lattice angles $\alpha_i$ for $i=1,\ldots, n$;\\
-- there exist a set of integers $M=\{m_1,\ldots, m_{n-1}\}$ such
that
$$
\sum\limits_{M,i=1}^{n}\bar{\pi{-}\alpha_i}= 2\pi.
$$
\end{theorem}

\begin{proof}
Consider an arbitrary $n$-tuple of ordinary lattice angles
$\alpha_i$, here $i=1,\ldots, n$.

Suppose that there exist a convex polygon $A_1A_2\ldots A_n$ with
consecutive angles $\alpha_i$ for $i=1,\ldots, n$. Let also the
couple of vectors $\bar{A_2A_3}$ and $\bar{A_2A_1}$ defines the
positive orientation of the plane (otherwise we apply to the
polygon $A_1A_2\ldots A_n$ some lattice-affine transformation
changing the orientation and come to the initial position).

Let $B_1=O{+}\bar{A_{n}A_{1}}$, and $B_i=O{+}\bar{A_{i-1}A_{i}}$
for $i=2,\ldots, n$. We put by definition
$$
\beta_i= \left\{
\begin{array}{ll}
\angle B_{i}OB_{i+1}, &\hbox{if $i=1,\ldots, n{-}1$}\\
\angle B_nOB_{1}, &\hbox{if $i=n$}\\
\end{array}
\right. .
$$
Consider the union of the sails for all $\beta_i$. This lattice
oriented broken line is of the class of the expanded lattice
angle with the normal form $2\pi{+}0$. The signed length-sine
sequence for this broken line contains exactly $n{-}1$ elements
that are not contained in the length-sine sequences for the sails
of $\beta_i$. Denote these numbers by $m_1, \ldots, m_{n-1}$, and
the set $\{m_1,\ldots, m_{n-1}\}$ by $M$. Then
$$
\sum\limits_{M,i=1}^{n}\bar{\beta_i}= 2\pi.
$$

From the definition of $\beta_i$ for $i=1,\ldots, n$ it follows
that $\beta_i \cong \pi {-}\alpha_{i}$. Therefore,
$$
\sum\limits_{M,i=1}^{n}\bar{\pi {-}\alpha_i}= 2\pi.
$$
We complete the proof of the statement in one side.

\vspace{2mm}

Suppose now, that there exist a set of integers $M=\{m_1,\ldots,
m_{n-1}\}$ such that
$$
\sum\limits_{M,i=1}^{n}\bar{\pi {-}\alpha_i}= 2\pi.
$$
This implies that there exist lattice points $B_1=(1,0)$,
$B_i=(x_i,y_i)$, for $i=2,\ldots {n{-}1}$, and $B_n=(-1,0)$ such
that
$$
\angle B_iOB_{i-1}\cong \pi {-}\alpha_{i-1} \hbox{, for
$i=2,\ldots, n$,} \quad \hbox{and}\quad \angle B_1OB_{n}\cong \pi
{-}\alpha_{n}.
$$

Denote by $M$ the lattice point
$$
O+\sum\limits_{i=1}^{n}\bar{OB_i}.
$$

Since all $\alpha_i$ are non-zero, the angles $ \pi {-}\alpha_i$
are ordinary. Hence, the origin $O$ is an interior point of the
convex hull of the points $B_i$ for $i=1,\ldots,k$. This implies
that there exist two consecutive lattice points $B_s$ and
$B_{s+1}$ (or $B_n$ and $B_1$), such that the lattice triangle
$\triangle B_s M B_{s+1}$ contains $O$ and the edge $B_s B_{s+1}$
does not contain $O$. Therefore,
$$
O=\lambda_1\bar{OM}+\lambda_2\bar{OB_i}+\lambda_3\bar{OB_{i+1}},
$$
where $\lambda_1$ is a positive integer, and $\lambda_2$ and
$\lambda_3$ are nonnegative integers. So there exist positive
integers $a_i$, where $i=1,\ldots, n$, such that
$$
O=O+\sum\limits_{i=1}^{n}(a_i\bar{OB_i}).
$$

Put by definition $A_0=O$, and $A_i=A_{i-1}+a_i\bar{OB_{i}}$ for
$i=2,\ldots, n$. The broken line $A_0A_1\ldots A_n$ is lattice
and by the above it is closed (i.~e. $A_0=A_n$). By construction,
the ordinary lattice angle at the vertex $A_i$ of the closed
lattice broken line is proper lattice-congruent to $\alpha_i$
($i=1,\ldots n$). Since the integers $a_i$ are positive for
$i=1,\ldots, n$ and the vectors $\bar{OB_i}$ are all in the
counterclockwise order, the broken line is a convex polygon.

The proof of Theorem~\ref{sum_for_polygons} is completed.
\end{proof}

\begin{remark}
Theorem~\ref{sum_for_polygons} generalizes the statement of
Theorem~\ref{sum}a. Note that the direct generalization of
Theorem~\ref{sum}b is false: the ordinary lattice angles do not
uniquely determine the proper lattice-affine homothety types of
convex polygons. See an example on Figure~\ref{not_b}.
\begin{figure}[h]
$$\epsfbox{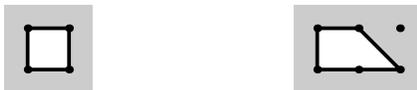}$$
\caption{An example of different types of polygons with the proper
lattice-congruent ordinary lattice angles.}\label{not_b}
\end{figure}
\end{remark}

\section{On lattice irrational case.}\label{iic}

The aim of this section is to generalize the notions of ordinary
and expanded lattice angles and their sums to the case of angles
with lattice vertices but not necessary lattice rays. We find
normal forms and extend the definition of lattice sums for a
certain special case of such angles.

\subsection{Infinite ordinary continued fractions.}

We start with the standard definition of infinite ordinary
continued fraction.
\begin{theorem}
Consider a sequence $(a_0, a_1, \ldots, a_n, \ldots )$ of
positive integers. There exists the following limit:
$r=\lim\limits_{k\to \infty}\bigl([a_0, a_1, \ldots ,a_k]\bigr)$.
\qed
\end{theorem}

This representation of $r$ is called an {\it$($infinite$)$
ordinary continued fraction} for $r$ and denoted by $[a_0, a_1,
\ldots ,a_n, \ldots ]$.

\begin{theorem}
For any irrational there exists and unique infinite ordinary
continued fraction. Any rational does not have infinite ordinary
continued fractions. \qed
\end{theorem}

For the proofs of these theorems we refer to the book~\cite{Hin}
by A.~Ya.~Hinchin.

\subsection{Ordinary lattice irrational angles.}

Let $A$, $B$, and $C$ do not lie in the same straight line.
Suppose also that $B$ is lattice. We denote the angle with the
vertex at $B$ and the rays $BA$ and $BC$ by $\angle ABC$. If the
open ray $BA$ contains lattice points, and the open ray $BC$ does
not contain lattice points, then we say that the angle $\angle
ABC$ is ordinary lattice {\it R-irrational} angle. If the open
ray $BA$ does not contain lattice points, and the open ray $BC$
contains lattice points, then we say that the angle $\angle AOB$
is ordinary lattice {\it L-irrational} angle. If the union of
open rays $BA$ and $BC$ does not contain lattice points, then we
say that the angle $\angle ABC$ is ordinary lattice {\it
LR-irrational} angle. We also call R-irrational, L-irrational,
and LR-irrational angles by {\it irrational} angles.

\begin{definition}
Two ordinary lattice irrational angles $\angle AOB$ and $\angle
A'O'B'$ are said to be {\it lattice-congruent} if there exist a
lattice-affine transformation which takes the vertex $O$ to the
vertex $O'$ and the rays $OA$ and $OB$ to the rays $O'A'$ and
$O'B'$ respectively. We denote this as follows: $\angle AOB \cong
\angle A'O'B'$.
\end{definition}

\subsection{Lattice-length sequences for ordinary lattice irrational angles.}

In this subsection we generalize the notion of lattice-length
sequences for the case of ordinary lattice irrational angles and
study its elementary properties.

Consider an ordinary lattice angle $\angle AOB$. Let also the
vectors $\bar{OA}$ and $\bar{OB}$ be linearly independent.

Denote the closed convex solid cone for the ordinary lattice
irrational angle $\angle AOB$ by $C(AOB)$. The boundary of the
convex hull of all lattice points of the cone $C(AOB)$ except the
origin is homeomorphic to the straight line. The closure in the
plane of the intersection of this boundary with the open cone
$AOB$ is called the {\it sail} for the cone $C(AOB)$. A lattice
point of the sail is said to be a {\it vertex} of the sail if
there is no lattice segment of the sail containing this point in
the interior. The sail of the cone $C(AOB)$ is a broken line with
an infinite number of vertices and without self-intersections. We
orient the sail in the direction from $\bar{OA}$ to $\bar{OB}$.
(For the definition of the sail and its higher dimensional
generalization, see, for instance, the works~\cite{Arn2},
\cite{Kor2}, and~\cite{Kar1}.)

In the case of ordinary lattice R-irrational angle we denote the
vertices of the sail by $V_i$, for $i\ge 0$, according to the
orientation of the sail (such that $V_0$ is contained in the ray
$OA$). In the case of ordinary lattice L-irrational angle we
denote the vertices of the sail by $V_{-i}$, for $i \ge 0$,
according to the orientation of the sail (such that $V_0$ is
contained in the ray $OB$). In the case of ordinary lattice
LR-irrational angle we denote the vertices of the sail by
$V_{-i}$, for $i \in \z$, according to the orientation of the
sail (such that $V_0$ is an arbitrary vertex of the sail).

\begin{definition}
Suppose that the vectors $\bar{OA}$ and $\bar{OB}$ of an ordinary
lattice irrational angle $\angle AOB$ are linearly independent.
Let $V_i$ be the vertices of the corresponding sail. The sequence
of lattice lengths and sines
$$
\begin{array}{l}
(\il(V_0V_1),\isin\angle V_0V_1V_2,\il(V_1V_2),\isin\angle V_1V_2V_3,\ldots), \hbox{ or}\\
(\ldots,\isin\angle V_{-3}V_{-2}V_{-1},\il(V_{-2}V_{-1}),\isin\angle V_{-2}V_{-1}V_{0},\il(V_{-1}V_{0})), \hbox{ or}\\
(\ldots,\isin\angle V_{-2}V_{-1}V_{0},\il(V_{-1}V_{0}),\isin\angle V_{-1}V_{0}V_{1},\il(V_{0}V_{1}),\ldots)\\
\end{array}
$$
is called the {lattice length-sine} sequence for the ordinary
lattice irrational angle $\angle AOB$, if this angle is
R-irrational, L-irrational, or LR-irrational respectively.
\end{definition}

\begin{proposition}
{\bf a).} The elements of the lattice length-sine sequence for
any ordinary
lattice irrational angle are positive integers.\\
{\bf b).} The lattice length-sine sequences of lattice-congruent
ordinary lattice irrational angles coincide. \qed
\end{proposition}

\subsection{Lattice tangents for ordinary lattice R-irrational angles.}

In this subsection we show, that the notion of lattice tangent is
well-defined for the case of ordinary lattice R-irrational angles.
We also formulate some basic properties of lattice tangent.

Let us generalize a notion of tangent to the case of R-irrational
angles using the property of Theorem~\ref{itan} for tangents of
ordinary angles.

\begin{definition}\label{itan_inf}
Let the vectors $\bar{OA}$ and $\bar{OB}$ of an ordinary lattice
R-irrational angle $\angle AOB$ be linearly independent. Suppose
that $V_i$ are the vertices of the corresponding sail. Let
$$
(\il(V_0V_1),\isin\angle V_0V_1V_2, \ldots, \isin\angle
V_{n-2}V_{n-1}V_n,\il(V_{n-1}V_n), \ldots)
$$
be the lattice length-sine sequence for the ordinary lattice
angle $\angle AOB$. The {\it lattice tangent} of the ordinary
lattice R-irrational angle $\angle AOB$ is the following number
$$
[\il(V_0V_1),\isin\angle V_0V_1V_2, \ldots, \isin\angle
V_{n-2}V_{n-1}V_n,\il(V_{n-1}V_n),\ldots].
$$
We denote it by $\itan\angle AOB$. We say also that this number is
the {\it continued fraction associated with the sail} of the
ordinary lattice R-irrational angle $\angle AOB$.
\end{definition}

\begin{proposition}
{\bf a).} For any ordinary lattice R-irrational angle $\angle
AOB$ with linearly independent
vectors $\bar{OA}$ and $\bar{OB}$ the number $\itan\angle AOB$ is irrational and greater than $1$.\\
{\bf b).} The values of the function $\itan$ at two
lattice-congruent ordinary lattice angles coincide. \qed
\end{proposition}

\subsection{Lattice arctangent for ordinary lattice R-irrational angles.}

We continue with the definition of lattice arctangents for
ordinary lattice R-irrational angles and the main properties of
these arctangents.

Consider the system of coordinates $OXY$ on the space $\r^2$ with
the coordinates $(x,y)$ and the origin $O$. We work with the
integer lattice of $OXY$.

For any rationals $p_1$ and $p_2$ we denote by $\alpha_{p_1,p_2}$
the angle with the vertex at the origin and two edges $\{(x,p_i
x)|x>0\}$, where $i=1,2$.

\begin{definition}\label{iarctan_inf}
For any irrational real $s> 1$, the ordinary lattice angle
$\angle AOB$ with the vertex $O$ at the origin, $A=(1,0)$, and
$B=(1,s)$, is called the {\it lattice arctangent} of $s$ and
denoted by $\iarctan s$.
\end{definition}

\begin{theorem}\label{itan_iarctan_inf}
{\bf a).} For any irrational $s$, such that $s> 1$,
$$
\itan(\iarctan s)=s.
$$
{\bf b).} For any ordinary lattice R-irrational angle $\alpha$
the following holds:
$$
\iarctan(\itan \alpha) \cong \alpha.
$$
\end{theorem}

\begin{proof}
Let us prove the first statement of the theorem. Let $s> 1$ be
some irrational real. Suppose that the sail of the angle
$\iarctan s$ is the infinite broken line $A_0A_1\ldots$ and the
corresponding ordinary continued fraction is
$[a_0,a_1,a_2,\ldots]$. Let also the coordinates of $A_i$ be
$(x_i,y_i)$.

We consider the ordinary lattice angles $\alpha_i$, corresponding
to the broken lines $A_0\ldots A_i$, for $i>0$. Then,
$$
\lim\limits_{i\to \infty} (y_i/x_i)=s/1.
$$

By Theorem~\ref{itan_iarctan} for any positive integer $i$ the
ordinary lattice angle $\alpha_i$ coincides with
$\iarctan([a_0,a_1,\ldots,a_{2i-2}])$, and hence the coordinates
$(x_i,y_i)$ of $A_i$ satisfy
$$
y_i/x_i=[a_0,a_1,\ldots,a_{2i-2}].
$$
Therefore,
$$
\lim\limits_{i\to \infty} ([a_0,a_1,\ldots,a_{2i-2}])=s.
$$
So, we obtain the first statement of the theorem:
$$
\itan(\iarctan s)=s.
$$

\vspace{2mm}

Now we prove the second statement. Consider an ordinary lattice
R-irrational angle $\alpha$. Suppose that the sail of the angle
$\alpha$ is the infinite broken line $A_0A_1\ldots$.

For any positive integer $i$ we consider the ordinary angle
$\alpha_i$, corresponding to the broken lines $A_0\ldots A_i$.

Denote by $C(\beta)$ the cone, corresponding to the ordinary
lattice (possible irrational) angle $\beta$. Note that
$C(\beta')$ and $C(\beta '')$ are lattice-congruent iff $\beta
\cong \beta'$.

By Theorem~\ref{itan_iarctan} we have:
$$
\iarctan(\itan \alpha_i) \cong \alpha_i.
$$
Since for any positive integer $n$ the following is true
$$
\bigcup\limits_{i=1}^{n}C(\alpha_i) \cong
\bigcup\limits_{i=1}^{n}C(\iarctan(\itan \alpha_i))
$$
we obtain
$$
C(\alpha)\cong\bigcup\limits_{i=1}^{\infty}C(\alpha_i) \cong
\bigcup\limits_{i=1}^{\infty}C(\iarctan(\itan \alpha_i)) \cong
C(\iarctan(\itan \alpha)).
$$
Therefore,
$$
\iarctan(\itan \alpha) \cong \alpha.
$$

This concludes the proof of Theorem~\ref{itan_iarctan_inf}.
\end{proof}

There is a description of ordinary lattice R-irrational angles
similar to the description of ordinary lattice angles (see
Proposition~\ref{basic_properties}e).

\begin{theorem}\label{angles_inf}{\bf (Description of ordinary lattice R-irrational angles.)}\\
{\bf a).} For any sequence of positive integers $(a_0, a_1, a_2,
\ldots)$ there exist
some ordinary  lattice R-irrational angle $\alpha$ such that $\itan\alpha=[a_0, a_1, a_2,\ldots ]$.\\
{\bf b).} Two ordinary lattice R-irrational angles are
lattice-congruent iff they have equal lattice tangents.
\end{theorem}

\begin{proof}
Theorem~\ref{itan_iarctan_inf}a implies the first statement of
the theorem.

Let us prove the second statement. Suppose that the ordinary
lattice R-irrational angles $\alpha$ and $\beta$ are
lattice-congruent, then their sails are also lattice-congruent.
Thus their lattice-sine sequences coincide. Therefore,
$\itan\alpha=\itan\beta$.

Suppose now that the lattice tangents for two ordinary lattice
R-irrational angles $\alpha$ and $\beta$ are equivalent. Since
for any irrational real the corresponding continued fraction
exists and is unique, the lattice length-sine sequences for the
the sails of $\alpha$ and $\beta$ coincide.

Let  $V_\alpha$ be the vertex of the angle $\alpha$, and
$A_0A_1\ldots$ be the sail for $\alpha$. For any positive integer
$i$ denote by $\alpha_i$ the angle $A_0V_{\alpha}A_i$ (with the
sail $A_0\ldots A_i$). Let  $V_\beta$ be the vertex of the angle
$\beta$, and $B_0B_1\ldots$ be the sail for $\beta$. For any
positive integer $i$ denote by $\beta_i$ the angle
$B_0V_{\alpha}B_i$ (with the sail $B_0\ldots B_i$).

Since the ordinary lattice R-irrational angles $\alpha$ and
$\beta$ have the same lattice length-sine sequences, for any
integer $i$ the ordinary lattice angles $\alpha_i$ and $\beta_i$
have the same lattice length-sine sequences, and, therefore,
$\alpha_i\cong\beta_i$.

Consider the lattice-affine transformation $\xi$ that takes the
vertex $V_\alpha$ to the vertex $V_\beta$, the lattice point
$A_0$ to the lattice point $B_0$ and $A_1$ to $B_1$. (Such
transformation exist since $\il(A_0A_1)=\il(B_0B_1)$.)

Choose an arbitrary $i\ge 1$. Since
$$
\xi(A_0)=B_0, \qquad \xi(A_1)=B_1, \quad  \hbox{and} \quad
\alpha_i\cong \beta_i,
$$
we have $\xi(A_i)=\xi(B_i)$. This implies that  the
lattice-affine transformation $\xi$ takes the sail for the
ordinary lattice R-irrational angle $\alpha$ to the sail for the
ordinary lattice R-irrational angle $\beta$. Therefore, the
angles $\alpha$ and $\beta$ are lattice-congruent.
\end{proof}

\begin{corollary}\label{angles_infL}{\bf (Description of ordinary lattice L-irrational
and LR-irrational angles.)}\\
{\bf a).} For any sequence of positive integers $(\ldots,
a_{-2},a_{-1},a_0)$ $($respectively $(\ldots a_{-1}, a_0, a_{1},
\ldots)$$)$ there exists an ordinary lattice L-irrational
$($LR-irrational$)$ angle with the given
LLS-sequence.\\
{\bf b).} Two ordinary lattice L-irrational  $($LR-irrational$)$
angles are lattice-congruent iff they have the same LLS-sequences.
\end{corollary}

\begin{proof}
The statement on L-irrational angles follows immediately from
Theorem~\ref{angles_inf}.

Let us construct a LR-angle with a given LLS-sequence $(\ldots
a_{-1}, a_0, a_{1}, \ldots)$. First we construct
$$
\alpha_1=\iarctan([a_0,a_1,a_2,\ldots]).
$$
Denote the points $(1,0)$ and $(1,a_0)$ by $A_0$ and $A_1$ and
construct the angle $\alpha_2$ lattice-congruent to the angle
$$
\iarctan([a_0,a_{-1},a_{-2},\ldots]).
$$
having the first two vertices $A_1$ and $A_0$ respectively. Now
the angle obtained by the rays of $\alpha_1$ and $\alpha_2$ that
do not contain lattice points is the LR-angle with the given
LLS-sequence.

Suppose now we have two LR-angles $\beta_1$ and $\beta_2$ with
the same LLS-sequences. Consider a lattice transformation taking
the vertex of $\beta_2$ to the vertex of $\beta_1$, and one of
the segments of $\beta_2$ to the segment $B_0B_1$ of $\beta_1$
with the corresponding order in LLS-sequence. Denote this angle
by $\beta_2'$. Consider the R-angles $\bar{\beta_1}$ and
$\bar{\beta_2'}$ corresponding to the sequences of vertices of
$\beta_1$ and $\beta_2'$ starting from $V_0$ in the direction to
$V_1$. These two angles are lattice-congruent by
Theorem~\ref{angles_inf}, therefore $\bar{\beta_1}$ and
$\bar{\beta_2'}$ coincide. So the angles $\beta_1$ and $\beta_2'$
have a common ray. By the same reason the second ray of $\beta_1$
coincides with the second ray of $\beta_2'$. Therefore $\alpha_1$
coincides with $\beta_2'$ and lattice-congruent to $\beta_2$.
\end{proof}

\subsection{Lattice signed length-sine infinite sequences.}

In this section we work in the oriented two-dimensional real
vector space with the fixed lattice. As previously, we fix
coordinates $OXY$ on this space.

A union of (ordered) lattice segments $\ldots, A_{i-1}A_i,
A_iA_{i+1}, A_{i+1}A_{i+2}, \ldots$ infinite to the right (to the
left, or both sides) is said to be a {\it lattice oriented
R-infinite $($L-infinite, or LR-infinite$)$ broken line}, if any
segment of the broken line is not of zero length, and any two
consecutive segments are not contained in the same straight line.
We denote this broken line by $\ldots A_{i-1}A_iA_{i+1}A_{i+2}
\ldots$ We also say that the lattice oriented broken line $\ldots
A_{i+2}A_{i+1}A_{i}A_{i-1} \ldots$ is {\it inverse} to the broken
line $\ldots A_{i-1}A_iA_{i+1}A_{i+2} \ldots$

\begin{definition}
Consider a lattice infinite oriented broken line and a point not
in this line. The broken line is said to be {\it on the unit
distance} from the point if all edges of the broken line are on
the unit lattice distance from the given point.
\end{definition}

Now, let us now associate to any lattice infinite oriented broken
line on the unit distance from some point the following sequence
of non-zero elements.

\begin{definition}
Let $\ldots A_{i-1}A_iA_{i+1}A_{i+2} \ldots$ be a lattice
oriented infinite broken line on the unit distance from some
lattice point $V$. Let
$$\begin{array}{l}
a_{2i-3}=\sgn(A_{i-2}VA_{i-1})\sgn(A_{i-1}VA_{i})\sgn(A_{i-2}A_{i-1}A_i)\isin\angle A_{i-2}A_{i-1}A_{i},\\
a_{2i-2}=\sgn(A_{i-1}VA_i)\il(A_{i-1}A_i)
\end{array}
$$
for all possible indexes $i$. The sequence $(\ldots a_{2i-3},
a_{2i-2}, a_{2i-1}\ldots)$ is called a {\it lattice signed
length-sine} sequence for the lattice oriented infinite broken
line on the unit distance from $V$.
\end{definition}

\begin{proposition}\label{ccc_inf}
A lattice signed length-sine sequence for the given lattice
infinite oriented broken line and the point is invariant under
the group action of proper lattice-affine transformations.
\end{proposition}

\begin{proof}
The statement of the proposition holds, since the functions
$\sgn$, $\il$, and $\isin$ are invariant under the group action
of proper lattice-affine transformations.
\end{proof}

\subsection{Proper lattice-congruence of lattice oriented infinite broken lines
on the unit distance from the lattice points.}

Let us formulate a necessary and sufficient conditions for two
lattice infinite oriented broken lines on the unit distance from
two lattice points to be proper lattice-congruent.

\begin{theorem}\label{twobrok_inf}
The lattice signed length-sine sequences of two lattice infinite
oriented broken lines on the unit distance from lattice points
$V_1$ and $V_2$ coincide, iff there exist proper lattice-affine
transformation taking the point $V_1$ to $V_2$ and one oriented
broken line to the other.
\end{theorem}

\begin{proof}
The lattice signed length-sine sequence for any lattice infinite
oriented broken line on the unit distance is uniquely defined,
and by Proposition~\ref{ccc} is invariant under the group action
of proper lattice-affine transformations. Therefore, the lattice
signed length-sine sequences for two proper lattice-congruent
broken lines coincide.

Suppose now that we have two lattice oriented infinite broken
lines $\ldots A_{i-1}A_iA_{i+1} \ldots$, and $\ldots
B_{i-1}B_iB_{i+1} \ldots$ on the unit distance from the points
$V_1$ and $V_2$, and with the same lattice signed length-sine
sequences. Consider the lattice-affine transformation $\xi$ that
takes the point $V_1$ to $V_2$, $A_i$ to $B_i$, and $A_{i+1}$ to
$B_{i+1}$ for some integer $i$. Since $\sgn (A_iVA_{i+1})=\sgn
(B_iVB_{i+1})$, the lattice-affine transformation $\xi$ is proper.
By Theorem~\ref{twobrok} the transformation $\xi$ takes any
finite oriented broken line $A_sA_{s+1} \ldots A_t$ containing the
segment $A_iA_{i+1}$ to the oriented broken line $B_sB_{s+1}
\ldots B_t$. Therefore, the transformation $\xi$ takes the
lattice oriented infinite broken lines $\ldots A_{i-1}A_iA_{i+1}
\ldots$ to the oriented broken line $\ldots B_{i-1}B_iB_{i+1}
\ldots$ and the lattice point $V_1$ to the lattice point $V_2$.

This concludes the proof of Theorem~\ref{twobrok_inf}.
\end{proof}

\subsection{Equivalence classes of almost positive lattice infinite oriented broken lines and corresponding
expanded lattice infinite angles.} We start this subsection with
the following general definition.

\begin{definition}\label{almost-positive}
We say that the lattice infinite oriented broken line on the unit
distance from some lattice point is {\it almost positive} if the
elements of the corresponding lattice signed length-sine sequence
are all positive, except for some finite number of elements.
\end{definition}

Let $l$ be the lattice (finite or infinite) oriented broken line
$\ldots A_{n-1} A_{n} \ldots A_m A_{m+1} \ldots$ Denote by
$l(-\infty,A_n)$ the broken line $\ldots A_{n-1} A_{n}$. Denote
by $l(A_m,+\infty)$ the broken line $A_m A_{m+1} \ldots$ Denote
by $l(A_n,A_m)$ the broken line $A_{n} \ldots A_m$.

\begin{definition}
Two lattice oriented infinite broken lines $l_1$ and $l_2$ on
unit distance from $V$ are said to be {\it equivalent} if there
exist two vertices $W_{11}$ and $W_{12}$ of the broken line $l_1$
and two vertices $W_{21}$ and $W_{22}$ of the broken line $l_2$
such that
the following three conditions are satisfied:\\
{\it i$)$} the broken line $l_1(-\infty,W_{11})$ coincides with the broken line $l_2(-\infty,W_{21})$;\\
{\it ii$)$} the broken line $l_1(W_{12},+\infty)$ coincides with the broken line $l_2(W_{22},+\infty)$;\\
{\it iii$)$} the closed broken line generated by
$l_1(W_{11},W_{12})$ and the inverse of $l_2(W_{21},W_{22})$ is
homotopy equivalent to the point on $\r^2\setminus \{V\}$.
\end{definition}

Now we give the definition of expanded lattice irrational angles.

\begin{definition}
An equivalence class of lattice R/L/LR-infinite oriented broken
lines on unit distance from $V$ containing the broken line $l$ is
called the {\it expanded lattice R/L/LR-infinite angle for the
equivalence class of $l$} at the {\it vertex} $V$ and denoted by
$\angle (V;l)$ (or, for short, {\it expanded lattice
R/L/LR-infinite angle}).
\end{definition}

\begin{definition}
Two expanded irrational lattice angles $\Phi_1$ and $\Phi_2$ are
said to be {\it proper lattice-congruent} iff there exist a
proper lattice-affine transformation sending the class of lattice
oriented broken lines corresponding to $\Phi_1$ to the class of
lattice oriented broken lines corresponding to $\Phi_2$. We
denote it by $\Phi_1\pcong \Phi_2$.
\end{definition}

\begin{remark}
Since all sails for ordinary lattice irrational angles are
lattice infinite oriented broken lines, the set of all ordinary
lattice irrational angles is naturally embedded into the set of
expanded lattice irrational angles. An ordinary lattice
irrational angle with a sail $S$ corresponds to the expanded
lattice irrational angle with the equivalence class of the broken
line $S$.
\end{remark}

\subsection{Revolution number for expanded lattice L- and R-irrational angles.}

Let us extend the revolution number to the case of almost
positive infinite oriented broken lines.

\begin{definition}
Let $\ldots A_{i-1}A_{i}A_{i+1}\ldots$ be some lattice R-, L- or
LR-infinite almost positive oriented broken line, and
$r=\{V+\lambda\bar{v}|\lambda\ge 0\}$ be the oriented ray for an
arbitrary vector $\bar v$ with the vertex at $V$. Suppose that
the ray $r$ does not contain the edges of the broken line, and
the broken line does not contain the vertex $V$. We call the
number
$$
\begin{array}{ll}
\lim\limits_{n\to +\infty}\# (r,V,A_0A_1\ldots A_n)& \hbox{if the broken line is R-infinite,}\\
\lim\limits_{n\to +\infty}\# (r,V,A_{-n}\ldots A_{-1}A_0)& \hbox{if the broken line is L-infinite,}\\
\lim\limits_{n\to +\infty}\# (r,V,A_{-n}A_{-n+1}\ldots  A_n)& \hbox{if the broken line is LR-infinite}\\
\end{array}
$$
the {\it intersection number} of the ray $r$ and the lattice
almost positive infinite oriented broken line broken line $\ldots
A_{i-1}A_{i}A_{i+1}\ldots$ and denote it by $\# (r,V, \ldots
A_{i-1}A_{i}A_{i+1}\ldots )$.
\end{definition}

\begin{proposition}\label{intersection_irr}
The intersection number of the ray $r$ and an almost positive
lattice infinite oriented broken line is well-defined.
\end{proposition}

\begin{proof}
Consider an almost positive lattice infinite oriented broken line
$l$. Let us show that the broken line $l$ intersects the ray $r$
only finitely many times.

By Definition~\ref{almost-positive} there exist vertices $W_1$
and $W_2$ of this broken line such that the signed lattice-sine
sequence for the lattice oriented broken line $l(-\infty,W_1)$
contains only positive elements, and the signed lattice-sine
sequence for the oriented broken line $l(W_2,+\infty)$ also
contains only positive elements.

The positivity of lattice-sine sequences implies that the lattice
oriented broken lines $l(-\infty,W_1)$, and $l(W_2,+\infty)$ are
the sails for some angles with the vertex $V$. Thus, these two
broken lines intersect the ray $r$ at most once each. Therefore,
the broken line $l$ intersects the ray $r$ at most once at the
part $l(-\infty,W_1)$, only a finite number times at the part
$l(W_1,W_2)$, and at most once at the part $l(W_2,+\infty)$.

So, the lattice infinite oriented broken line $l$ intersects the
ray $r$ only finitely many times, and, therefore, the
corresponding intersection number is well-defined.
\end{proof}

Now we give a definition of the lattice revolution number for
expanded lattice R-irrational and L-irrational angles.

\begin{definition}
{\bf a).} Consider an arbitrary R-infinite (or L-infinite)
expanded lattice angle $\angle (V,l)$, where $V$ is some lattice
point, and $l$ is a lattice infinite oriented almost-positive
broken line. Let $A_0$ be the first (the last) vertex of $l$.
Denote the rays $\{V+\lambda\bar{VA_0}|\lambda\ge 0\}$ and
$\{V-\lambda\bar{VA_0}|\lambda\ge 0\}$ by $r_+$ and $r_-$
respectively. The following number
$$
\frac{1}{2}\bigl(\#(r_+,V,l)+\#(r_-,V,l)\bigr)
$$
is called the {\it lattice revolution number} for the expanded
lattice irrational angle $\angle (V,l)$, and denoted by
$\#(\angle(V,l)).$
\end{definition}

\begin{proposition}
The revolution number of an R-irrational $($or L-irrational$)$
expanded lattice angle is well-defined.
\end{proposition}

\begin{proof}
Consider an arbitrary expanded lattice R-irrational angle $\angle
(V, A_0A_1\ldots)$. Let
$$
r_+=\{V+\lambda\bar{VA_0}|\lambda\ge 0\} \quad  \hbox{and} \quad
r_-=\{V-\lambda\bar{VA_0}|\lambda\ge 0\}.
$$

Since the lattice oriented broken line $A_0A_1A_2\ldots$ is on
the unit lattice distance from the point $V$, any segment of this
broken line is on the unit lattice distance from $V$. Thus, the
broken line does not contain $V$, and the rays $r_+$ and $r_-$ do
not contain edges of the broken line.

Suppose that
$$
\angle V,A_0A_1A_2\ldots =\angle V',A'_0A_1'A_2'\ldots
$$
This implies that $V=V'$, $A_0=A'_0$, $A_{n+k}=A_{m+k}'$ for some
integers $n$ and $m$ and any non-negative integer $k$, and the
broken lines $A_0A_1\ldots A_nA_{m-1}'\ldots A_1' A_0'$ is
homotopy equivalent to the point on $\r^2\setminus \{V\}$.  Thus,
$$
\begin{array}{l}
\#(\angle V,A_0A_1\ldots)-\#(\angle V',A'_0A_1'\ldots)=\\
\frac{1}{2}\bigl(\#(r_+,A_0A_1\ldots A_nA_{m-1}'\ldots A_1' A_0')+
\#(r_-,A_0A_1\ldots A_nA_{m-1}'\ldots A_1' A_0')\bigr)\\
=0{+}0=0.
\end{array}
$$
And hence
$$\
\#(\angle V,A_0A_1A_2\ldots)=\#(\angle V',A'_0A_1'A_2'\ldots ).
$$

Therefore, the revolution number of any expanded lattice
R-irrational angle is well-defined.

The proof for L-irrational angles repeats the proof for
R-irrational angles and is omitted here.
\end{proof}

\begin{proposition}
The revolution number of expanded lattice R/L-irrational angles
is invariant under the group action of the proper lattice-affine
transformations. \qed
\end{proposition}

\subsection{Normal forms of expanded lattice R- and L-irrational lattice angles.}

In this subsection we formulate and prove a theorem on normal
forms of expanded lattice R-irrational and L-irrational lattice
angles.

For the theorems of this subsection we introduce the following
notation. By the sequence
$$
\bigl((a_0,\ldots,a_n)\times \hbox{$k$-times},b_0,b_1\ldots\bigr),
$$
where $k\ge 0$, we denote the sequence:
$$
(\underbrace{a_0,\ldots,a_n,a_0,\ldots,a_n, \quad \ldots\quad
,a_0,\ldots,a_n}_{\mbox{$k$-times}} ,b_0,b_1,\ldots).
$$
By the sequence
$$
\bigl(\ldots,b_{-2},b_{-1},b_0,(a_0,\ldots,a_n)\times
\hbox{$k$-times}\bigr),
$$
where $k\ge 0$, we denote the following sequence:
$$
(\ldots,b_{-2},b_{-1},b_0,\underbrace{a_0,\ldots,a_n,a_0,\ldots,a_n,
\quad \ldots\quad ,a_0,\ldots,a_n}_{\mbox{$k$-times}}).
$$

We start with the case of expanded lattice R-irrational angles.

\begin{definition}\label{normal_forms_r_inf}
Consider a lattice R-infinite oriented broken line $A_0A_1\ldots$
on the unit distance from the origin $O$. Let also $A_0$ be the
point $(1,0)$, and the point $A_1$ be on the line $x=1$. If the
signed length-sine sequence of the expanded ordinary R-irrational
angle $\Phi_0=\angle(O,A_0A_1\ldots )$
coincides with the following sequence (we call it {\it characteristic sequence} for the corresponding angle):\\
{\bf IV$_k$)} $\bigl((1,-2,1,-2)\times \hbox{$k$-times},
a_0,a_1,\ldots \bigr)$, where $k\ge 0$, $a_i>0$, for $i\ge 0$,
then we denote the angle $\Phi_0$ by $k\pi {+}
\iarctan([a_0,a_1,\ldots])$
and say that $\Phi_0$ is {\it of the type} {\bf IV$_k$};\\
{\bf V$_k$)} $\bigl((-1,2,-1,2)\times \hbox{$k$-times},
a_0,a_1,\ldots\bigr)$, where $k> 0$, $a_i> 0$, for $i\ge 0$, then
we denote the angle $\Phi_0$ by $-k\pi +
\iarctan([a_0,a_1,\ldots])$ and say that $\Phi_0$ is {\it of the
type} {\bf V$_k$}.
\end{definition}

\begin{theorem}\label{nf_r_inf}
For any expanded lattice R-irrational angle $\Phi$ there exist
and unique a type among the types {\bf{IV}}-{\bf{V}} and a unique
expanded lattice R-irrational angle $\Phi$ of that type such that
$\Phi_0$ is proper lattice
congruent to $\Phi_0$.\\
{\rm The expanded lattice R-irrational angle $\Phi_0$ is said to
be {\it the normal form} for the expanded lattice R-irrational
angle $\Phi$.}
\end{theorem}

\begin{proof}
First, we prove that any two distinct expanded lattice
R-irrational angles listed in Definition~\ref{normal_forms_r_inf}
are not proper lattice-congruent. Let us note that the revolution
numbers of expanded lattice angles distinguish the types of the
angles. The revolution number for the expanded lattice angles of
the type {\bf IV$_k$} is $1/4{+}1/2k$ where $k\ge 0$. The
revolution number for the expanded lattice angles of the type
{\bf V$_k$} is $1/4{-}1/2k$ where $k > 0$.

We now prove that the normal forms of the same type {\bf IV$_k$}
(or {\bf V$_k$}) are not proper lattice-congruent for any integer
$k$. Consider the expanded lattice R-infinite angle $\Phi= k\pi
{+} \iarctan([a_0,a_1,\ldots])$. Suppose that a lattice oriented
broken line $A_0A_1A_2 \ldots$ on the unit distance from $O$
defines the angle $\Phi$. Let also that the signed lattice-sine
sequence for this broken line be characteristic.

If $k$ is even, then the ordinary lattice R-irrational angle with
the sail $A_{2k}A_{2k+1}\ldots$ is proper lattice-congruent to
the angle $\iarctan([a_0,a_1,\ldots])$. This angle is a proper
lattice-affine invariant for the expanded lattice R-irrational
angle $\Phi$ (since $A_{2k}=A_0$). This invariant distinguish the
expanded lattice R-irrational angles of type {\bf IV$_k$} (or
{\bf V$_k$}) for even $k$.

If $k$ is odd, then denote $B_i=V{+}\bar{A_iV}$. The ordinary
lattice R-irrational angle with the sail $B_{2k}B_{2k+1}\ldots$
is proper lattice-congruent to the angle
$\iarctan([a_0,a_1,\ldots])$. This angle is a proper
lattice-affine invariant of the expanded lattice R-irrational
angle $\Phi$ (since $B_{2k}=V{+}\bar{A_0V}$). This invariant
distinguish the expanded lattice R-irrational angles {\bf IV$_k$}
(or {\bf V$_k$}) for odd $k$.

Therefore, the expanded lattice angles listed in
Definition~\ref{normal_forms_r_inf} are not proper
lattice-congruent.

\vspace{2mm}

Secondly, we prove that an arbitrary expanded lattice
R-irrational angle is proper lattice-congruent to some of the
expanded lattice angles listed in
Definition~\ref{normal_forms_r_inf}.

Consider an arbitrary expanded lattice R-irrational angle
$\Phi=\angle (V,A_0A_1\ldots )$. Suppose that
$\#(\Phi)=1/4{+}k/2$ for some non-negative integer $k$. By
Proposition~\ref{intersection_irr} there exist an integer
positive number $n_0$ such that the lattice oriented broken line
$A_{n_0}A_{n_0+1}\ldots$ does not intersect the rays
$r_+=\{V{+}\lambda\bar{VA_0}|\lambda\ge 0\}$ and
$r_-=\{V{-}\lambda\bar{VA_O}|\lambda\ge 0\}$, and the signed
lattice length-sine sequence $(a_{2n_0-2},a_{2n_0-1}\ldots)$ for
the oriented broken line $A_{n_0}A_{n_0+1}\ldots$ does not
contain non-positive elements.

By Theorem~\ref{nf} there exist integers $k$ and $m$, and a
lattice oriented broken line
$$
A_0B_1B_2\ldots B_{2k}B_{2k+1}\ldots B_{2k+m}A_{n_0}
$$
with lattice length-sine sequence of the form
$$
\bigl((1,-2,1,-2)\times \hbox{$k$-times}, b_0,b_1, \ldots,
b_{2m-2} \bigr),
$$
where all $b_i$ are positives.

Consider now the lattice oriented infinite broken line
$A_0B_1B_2\ldots B_{2k+m-1}A_{n_0}A_{n_0+1}\ldots$ The
length-sine sequence for this broken line is as follows
$$
\bigl((1,-2,1,-2)\times \hbox{$k$-times}, b_0,b_1, \ldots,
b_{2m-2},v ,a_{2n_0-2}a_{2n_0-1}\ldots\bigr),
$$
where $v$ is (not necessary positive) integer.

Note that the lattice oriented broken line $A_0B_1B_2\ldots
B_{2k+m}A_{n_0}$ is a sail for the angle $\angle A_0VA_{n_0}$ and
the sequence $A_{n_0}A_{n_0+1}\ldots$ is a sail for some
R-irrational angle (we denote it by $\alpha$). Let $H_1$ be the
convex hull of all lattice points of the angle $\angle
A_0VA_{n_0}$ except the origin, and $H_2$ be the convex hull of
all lattice points of the angle $\alpha$ except the origin. Note
that $H_1$ intersects $H_2$ in the ray with the vertex at
$A_{n_0}$.

The lattice oriented infinite broken line $B_{2k}B_{2k+2}\ldots
B_{2k+m}A_{n_0}A_{n_0+1}\ldots$ intersects the ray $r_+$ in the
unique point $B_{2k}$ and does not intersect the ray $r_-$. Hence
there exists a straight line $l$ intersecting both boundaries of
$H_1$ and $H_2$, such that the open half-plane with the boundary
straight line $l$ containing the origin does not intersect the
sets $H_1$ and $H_2$.

Denote $B_0=A_0$ and $B_{2k+m+1}=A_{n_0}$. The intersection of
the straight line $l$ with $H_1$ is either a point $B_s$ (for $2k
\le s \le 2k+m+1$), or a boundary segment $B_sB_{s+1}$ for some
integer $s$ satisfying $2k\le s\le 2k{+}m$. The intersection of
$l$ with $H_2$ is either a point $A_t$ for some integer $t\ge
n_0$, or a boundary segment $A_{t-1}A_{t}$ for some integer
$t>n_0$.

Since the triangle $\triangle VA_tB_s$ does not contain interior
point of $H_1$ and $H_2$, the lattice points of $\triangle
VA_tB_s$ distinct to $B$ are on the segment $A_tB_s$. Hence, the
segment $A_tB_s$ is on unit lattice distance to the vertex $V$.
Therefore, the lattice infinite oriented broken line
$$
A_0B_1B_2\ldots B_sA_t A_{t+1}\ldots
$$
is on lattice unit distance.

Since the lattice oriented broken line $B_k\ldots B_{s}A_t
A_{t+1}\ldots$ is convex, it is a sail for some lattice
R-irrational angle. (Actually, the case $B_s=A_t=A_{n_0}$ is also
possible, then delete one of the copies of $A_{n_0}$ from the
sequence.) We denote this broken like by $C_{2k+1}C_{2k+2}\ldots$
The corresponding signed lattice length-sine sequence is
$(c_{4k},c_{4k+1},c_{4k+2},\ldots)$, where $c_i>0$ for $i \ge 4k$.
Thus the signed lattice length-sine sequence for the lattice
ordered broken line $A_0B_1B_2\ldots B_{2k}C_{2k+1}C_{2k+2}\ldots$
is
$$
\bigl((1,-2,1,-2)\times
\hbox{$(k-1)$-times},1,-2,1,w,(c_{4k},c_{4k+1},c_{4k+2},\ldots),
$$
where $w$ is an integer that is not necessary equivalent to $-2$.

Consider an expanded lattice angle $\angle (V,A_0B_1B_2\ldots
B_{2k}C_{2k+1})$. By Lemma~\ref{lema_nf} there exists a lattice
oriented broken line $C_0 \ldots C_{2k+1}$ with the vertices
$C_0=A_0$ and $C_{2k+1}$ of the same equivalence class, such that
$C_{2k}=B_{2k}$, and the signed lattice length-sine sequence for
it is
$$
\bigl((1,-2,1,-2)\times \hbox{$k$-times},c_{4k},c_{4k+1}).
$$
Therefore, the lattice oriented R-infinite broken line
$C_0C_1\ldots$ for the angle $\angle (V,A_0A_1\ldots)$ has the
signed lattice length-sign sequence coinciding with the
characteristic sequence for the angle
$k\pi+\iarctan([c_{4k},c_{4k+1},\ldots])$. Therefore,
$$
\Phi\pcong k\pi+\iarctan([c_{4k},c_{4k+1},\ldots]).
$$

This concludes the proof of the theorem for the case of
nonnegative integer $k$.

The proof for the case of negative $k$ repeats the proof for the
nonnegative case and is omitted here.
\end{proof}

Let us give the definition of trigonometric functions for
expanded lattice R-irrational angles.

\begin{definition}
Consider an arbitrary expanded lattice R-irrational angle $\Phi$
with
the normal form $k\pi {+} \varphi$ for some integer $k$.\\
{\bf a).} The ordinary lattice R-irrational angle $\varphi$ is
said to be
{\it associated} with the expanded lattice R-irrational angle $\Phi$.\\
{\bf b).} The number $\itan(\varphi)$ is called the lattice {\it
tangent} of the expanded lattice R-irrational angle $\Phi$.
\end{definition}

\vspace{2mm}

We continue now with the case of expanded lattice L-irrational
angles.

\begin{definition}
The expanded lattice irrational angle $\angle (V,\ldots A_{i+2}
A_{i+1} A_i \ldots)$ is said to be transpose to the expanded
lattice irrational angle $\angle (V,\ldots A_i A_{i+1} A_{i+2}
\ldots)$ and denoted by $\left( \angle (V,\ldots A_i A_{i+1}
A_{i+2} \ldots)\right)^t$.
\end{definition}

\begin{definition}\label{normal_forms_l_inf}
Consider a lattice L-infinite oriented broken line $\ldots
A_{-1}A_0$ on the unit distance from the origin $O$. Let also
$A_0$ be the point $(1,0)$, and the point $A_{-1}$ be on the
straight line $x=1$. If the signed length-sine sequence of the
expanded ordinary L-irrational angle $\Phi_0=\angle(O,\ldots
A_{-1}A_0 )$
coincides with the following sequence (we call it {\it characteristic  sequence} for the corresponding angle):\\
{\bf IV$_k$)} $\bigl(\ldots,a_{-1},a_{0},(-2,1,-2,1)\times
\hbox{$k$-times} \bigr)$, where $k\ge 0$, $a_i>0$, for $i\le 0$,
then we denote the angle $\Phi_0$ by $k\pi
{+}\iarctan^t([a_{0},a_{-1},\ldots])$
and say that $\Phi_0$ is {\it of the type} {\bf IV$_k$};\\
{\bf V$_k$)} $\bigl(\ldots,a_{-1},a_{0},(2,-1,2,-1)\times
\hbox{$k$-times}\bigr)$, where $k> 0$, $a_i>0$, for $i\le 0$, then
we denote the angle $\Phi_0$ by $-k\pi {+}
\iarctan^t([a_{0},a_{-1},\ldots])$ and say that $\Phi_0$ is {\it
of the type} {\bf V$_k$}.
\end{definition}

\begin{theorem}\label{nf_l_inf}
For any expanded lattice L-irrational angle $\Phi$ there exist
and unique a type among the types {\bf{IV}}-{\bf{V}} and a unique
expanded lattice L-irrational angle $\Phi_0$ of that type such
that $\Phi$ is proper lattice
congruent to $\Phi_0$.\\
{\rm The expanded lattice L-irrational angle $\Phi_0$ is said to
be {\it the normal form} for the expanded lattice L-irrational
angle $\Phi$.}
\end{theorem}

\begin{proof}
After transposing the set of all angles and change of the
orientation of the plane the statement of Theorem~\ref{nf_l_inf}
coincide with the statement of Theorem~\ref{nf_r_inf}.
\end{proof}

\subsection{Sums of expanded lattice angles and expanded lattice irrational angles.}

We conclude this section with a particular definitions of sums of
ordinary lattice angles, and ordinary lattice R-irrational or/and
L-irrational angles.

\begin{definition}
Consider expanded lattice angles $\Phi_i$, where $i=1,\ldots,t$,
an expanded lattice R-irrational angle $\Phi_r$, and an expanded
lattice L-irrational angle $\Phi_l$. Let the characteristic
signed lattice lengths-sine sequences for the normal forms of the
angles $\Phi_i$ be $(a_{0,i},a_{1,i},\ldots, a_{2n_i,i})$; of
$\Phi_r$ be $(a_{0,r},a_{1,r},\ldots)$, and of
$\Phi_l$ be $(\ldots,a_{-1,l},a_{0,l})$.\\
Let $M_R=(m_1,\ldots,m_{t-1},m_r)$ be some $t$-tuple of integers.
The normal form of any expanded lattice angle, corresponding to
the following lattice signed length-sine sequence
$$
\begin{aligned}
\bigl(a_{0,1},a_{1,1},\ldots, a_{2n_1,1},m_1,
a_{0,2},a_{1,2},\ldots, a_{2n_2,2},m_2, \ldots \\
\ldots ,m_{t-1},a_{0,t},a_{1,t},\ldots, a_{2n_t,t}
m_r,a_{0,r},a_{1,r},\ldots\bigr)
\end{aligned}
$$
is called the {\it $M_R$-sum of expanded lattice angles $\Phi_i$ $($$i=1,\ldots,t$$)$} and $\Phi_r$.\\
Let $M_L=(m_l,m_1,\ldots,m_{t-1})$ be some $t$-tuple of integers.
The normal form for any expanded lattice angle, corresponding to
the following lattice signed length-sine sequence
$$
\begin{aligned}
\bigl(\ldots,a_{-1,l},a_{0,l},m_l,a_{0,1},a_{1,1},\ldots,
a_{2n_1,1},m_1,
a_{0,2},a_{1,2},\ldots, a_{2n_2,2},m_2, \ldots \\
\ldots ,m_{t-1},a_{0,t},a_{1,t},\ldots, a_{2n_t,t}\bigr)
\end{aligned}
$$
is called the {\it $M_L$-sum of expanded lattice angles $\Phi_l$, and $\Phi_i$ $($$i=1,\ldots,t$$)$}.\\
Let $M_{LR}=(m_l,m_1,\ldots,m_{t-1},m_r)$ be some $(t+1)$-tuple
of integers. Any expanded lattice LR-irrational angle,
corresponding to the following lattice signed length-sine sequence
$$
\begin{aligned}
\bigl(\ldots,a_{-1,l},a_{0,l},m_l,a_{0,1},a_{1,1},\ldots,
a_{2n_1,1},m_1,
a_{0,2},a_{1,2},\ldots, a_{2n_2,2},m_2, \ldots \\
\ldots ,m_{t-1},a_{0,t},a_{1,t},\ldots, a_{2n_t,t}
m_r,a_{0,r},a_{1,r},\ldots\bigr)
\end{aligned}
$$
is called a {\it $M_{LR}$-sum of expanded lattice angles
$\Phi_l$, $\Phi_i$ $($$i=1,\ldots,t$$)$} and $\Phi_r$.
\end{definition}

\appendix

\section{On global relations on algebraic singularities of complex projective toric varieties
corresponding to integer-lattice triangles.}\label{toric}

In this appendix we describe an application of theorems on sums
of lattice tangents for the angles of lattice triangles and
lattice convex polygons to theory of complex projective toric
varieties. We refer the reader to the general definitions of
theory of toric varieties to the works of
V.~I.~Danilov~\cite{Dan}, G.~Ewald~\cite{Ewa},
W.~Fulton~\cite{Ful}, and T.~Oda~\cite{Oda}.

Let us briefly recall the definition of complex projective toric
varieties associated to lattice convex polygons. Consider a
lattice convex polygon $P$ with vertices $A_0, A_1,\ldots ,A_n$.
Let the intersection of this (closed) polygon with the lattice
consists of the points $B_i=(x_i,y_i)$ for $i=0, \ldots, m$. Let
also $B_i=A_i$ for $i=0, \ldots, n$. Denote by $\Omega$ the
following set in $\c P^{m}$:
$$
\Bigl\{ \bigl(t_1^{x_1}t_2^{y_1}t_3^{-x_1-y_1}:
      t_1^{x_2}t_2^{y_2}t_3^{-x_2-y_2}: \ldots:
      t_1^{x_m}t_2^{y_m}t_3^{-x_m-y_m}\bigr)| t_1,t_2,t_3 \in \c\setminus \{0\}
\Bigr\}.
$$
The closure of the set $\Omega$ in the natural topology of $\c
P^{m}$ is called the {\it complex toric variety associated with
the polygon $P$} and denoted by $X_P$.

For any $i=0,\ldots,m$ we denote by $\tilde A_i$ the point
$(0:\ldots:0:1:0:\ldots:0)$ where $1$ stands on the $(i{+}1)$-th
place.

From general theory it follows that:

$\bf a)$ the set $X_P$ is a complex projective
complex-two-dimensional variety with isolated algebraic
singularities;

$\bf b)$  the complex toric projective variety contains the
points $\tilde A_i$ for $i=0,\ldots,n$ (where $n{+}1$ is the
number of vertices of convex polygon);

$\bf c)$  the points of $X_P\setminus \{\tilde A_0, \tilde
A_1,\ldots ,\tilde A_n\}$ are non-singular;

$\bf d)$  the point $\tilde A_i$ for any integer $i$ satisfying
$0\le i\le n$ is singular iff the corresponding ordinary lattice
angle $\alpha_i$ at the vertex $A_i$ of the polygon $P$ is not
lattice-congruent to $\iarctan (1)$;

$\bf e)$ the algebraic singularity at $\tilde A_i$ for any
integer $i$ satisfying $0\le i\le n$ is uniquely determined by
the lattice-affine type of the non-oriented sail of the lattice
angle $\alpha_i$.

\vspace{2mm}

The algebraic singularity is said to be {\it toric} if there
exists a projective toric variety with the given algebraic
singularity.

Note that the lattice-affine classes of non-oriented sails for
angles $\alpha$ and $\beta$ coincide iff $\beta \cong \alpha$, or
$\beta \cong \alpha ^t$. This allows us to associate to any
complex-two-dimensional toric algebraic singularity,
corresponding to the sail of the angle $\alpha$, the unordered
couple of rationals $(a,b)$, where $a=\itan\alpha$ and
$b=\itan\alpha^t$.

\begin{remark}
Note that the continued fraction for the sail $\alpha$ is
slightly different to the Hirzebruch-Jung continued fractions for
toric singularities (see the works~\cite{Jun} by H.W.E.~Jung,
and~\cite{Hir} by F.~Hirzebruch). The relations between these
continued fractions is described in the paper~\cite{P-P} by
P.Popescu-Pampu.
\end{remark}

\begin{corollary}
Suppose, that we are given by three complex-two-dimensional
toric singularities defined by couples of rationals $(a_i,b_i)$
for $i=1,2,3$. There exist a complex toric variety associated
with some triangle with these three singularities iff there exist
a permutation $\sigma\in S_3$ and the rationals $c_i$ from the
sets $\{a_i,b_i\}$ for $i=1,2,3$, such that the following
conditions hold:

{\it i$)$} the rational $]c_{\sigma(1)} ,-1,c_{\sigma(2)}[$ is
either negative or greater than $c_{\sigma(1)}$;

{\it ii$)$} $]c_{\sigma(1)},-1,c_{\sigma(2)},-1,c_{\sigma(3)}[=0.$
\end{corollary}

We note again that we use odd continued fractions for $c_1$,
$c_2$, and $c_3$ in the statement of the above proposition (see
Subsection~\ref{q1q2q3} for the notation of continued fractions).

\begin{proof}
The proposition follows directly from Theorem~\ref{sum}a.
\end{proof}

\begin{proposition}\label{polygon_toric}
For any collection $($with multiplicities$)$ of
complex-two-dimensional toric algebraic singularities there exist
a complex-two-dimensional toric projective variety with exactly
the given collection of toric singularities.
\end{proposition}

For the proof of Proposition~\ref{polygon_toric} we need the
following lemma.

\begin{lemma}\label{polygon_toric_lemma}
For any collection of ordinary lattice angles $\alpha_i$
$($$i=1,\ldots, n$$)$, there exist an integer $k\ge n{-}1$ and a
k-tuple of integers $M=(m_1,\ldots, m_k)$, such that
$$
\bar{\alpha_1}+_{m_{1}}\ldots+_{m_{n-1}}\bar{\alpha_n}+_{m_n}\iarctan(1)+_{m_{n+1}}\ldots
+_{m_k}\iarctan(1) = 2\pi.
$$
\end{lemma}

\begin{proof}
Consider any collection of ordinary lattice angles $\alpha_i$
($i=1,\ldots, n$) and denote
$$
\Phi=\bar{\alpha_1}+_{1}\bar{\alpha_2}+_{1}\ldots
+_{1}\bar{\alpha_n}.
$$
There exist an oriented lattice broken line for the angle $\Phi$
with the signed lattice-signed sequence with positive elements.
Hence, $\Phi \pcong \varphi + 0\pi$.

If $\varphi\cong \iarctan(1)$, we have
$$
\Phi+_{-2}\iarctan(1)+_{-2}\iarctan(1)+_{-2}\iarctan(1)= 2\pi.
$$
Then $k=n+2$, and $M=(1,\ldots, 1,-2,-2,-2)$.

Suppose now $\varphi \not \cong \iarctan(1)$, then the following
holds
$$
\bar{\varphi}+_{-1}\bar{\pi
{-}\varphi}+_{-2}\iarctan(1)+_{-2}\iarctan(1)= 2\pi.
$$
Consider the sail for the angle $\pi{-}\varphi$. Suppose the
sequence of all its lattice points (not only vertices) is $B_0,
\ldots, B_s$ (with the order coinciding with the order of the
sail). Then we have
$$
\angle B_iOB_{i+1}\pcong \iarctan (1) \quad \hbox{for any
$i=1,\ldots,s$}.
$$
Denote by $b_i$ the values of $\isin\angle B_iOB_{i+1}$ for
$i=1,\ldots,s$. Then we have
$$
\begin{array}{l}
\bar{\varphi}+_{-2}\iarctan(1)+_{-2}\iarctan(1)+_{-2}\iarctan(1)=\\
\bar{\alpha_1}+_{1}\bar{\alpha_2}+_{1}\ldots +_{1}\bar{\alpha_n}+_{-1}\\
\iarctan(1)+_{b_1}\iarctan(1)+_{b_2}\ldots+_{b_{s}}\iarctan(1)+_{-2}\\
\iarctan(1)+_{-2}\iarctan(1)+_{-2}\iarctan(1) = 2\pi.
\end{array}
$$
Therefore, $k=n{+}s{+}3$, and
$$
M=(\underbrace{1,1,\ldots,1,1}_{\mbox{$(n{-}1)$-times}},-1,b_1,\ldots,b_s,-2,-2,-2).
$$

The proof of Lemma~\ref{polygon_toric_lemma} is completed.
\end{proof}

{\it Proof of the statement of the
Proposition~\ref{polygon_toric}.} Consider an arbitrary
collection of two-dimensional toric algebraic singularities.
Suppose that they are represented by ordinary lattice angles
$\alpha_i$ ($i=1,\ldots, n$). By Lemma~\ref{polygon_toric_lemma}
there exist an integer $k\ge n{-}1$ and a k-tuple of integers
$M=(m_1,\ldots, m_k)$, such that
$$
\bar{(\pi-\alpha_1)}+_{m_{1}}\ldots+_{m_{n-1}}\bar{(\pi-\alpha_n)}+_{m_{n}}
\iarctan(1)+_{m_{n+1}}\ldots+_{m_k}\iarctan(1) = 2\pi.
$$

By Theorem~\ref{sum_for_polygons} there exist a convex polygon
$P=A_0\ldots A_k$ with angles proper lattice-congruent to the
ordinary lattice angles $\alpha_i$ ($i=1,\ldots, n$), and
$k{-}n{+}1$ angles $\iarctan (1)$.

By the above, the toric variety $X_P$ is nonsingular at points of
$P_X\setminus \{\tilde A_0, \tilde A_1,\ldots ,\tilde A_k\}$. It
is also nonsingular at the points $\tilde A_i$ with the
corresponding ordinary lattice angles lattice-congruent to
$\iarctan(1)$. The collection of the toric singularities at the
remaining points coincide with the given collection.

This concludes the proof of Proposition~\ref{polygon_toric}. \qed

\begin{figure}[h]
$$\epsfbox{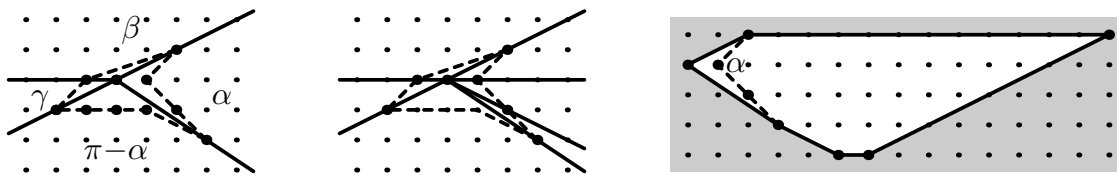}$$
\caption{Constructing a polygon with all angles proper
inte\-ger-con\-gruent to $\iarctan (1)$ except one angle that is
proper inte\-ger-con\-gruent to $\iarctan
(7/5)$.}\label{mnogoug_7/5}
\end{figure}

On Figure~\ref{mnogoug_7/5} we show an example of the polygon for
a projective toric variety with the unique toric singularity,
represented by the sail of  $\iarctan (7/5)$. The ordinary
lattice angle $\alpha$ on the figure is proper lattice-congruent
to $\iarctan (7/5)$, the angles $\beta$ and $\gamma$ are proper
lattice-congruent to $\iarctan (1)$.

\section{On lattice-congruence criterions for lattice triangles. Examples of lattice triangles.}\label{criterio}

\begin{figure}
$$\epsfbox{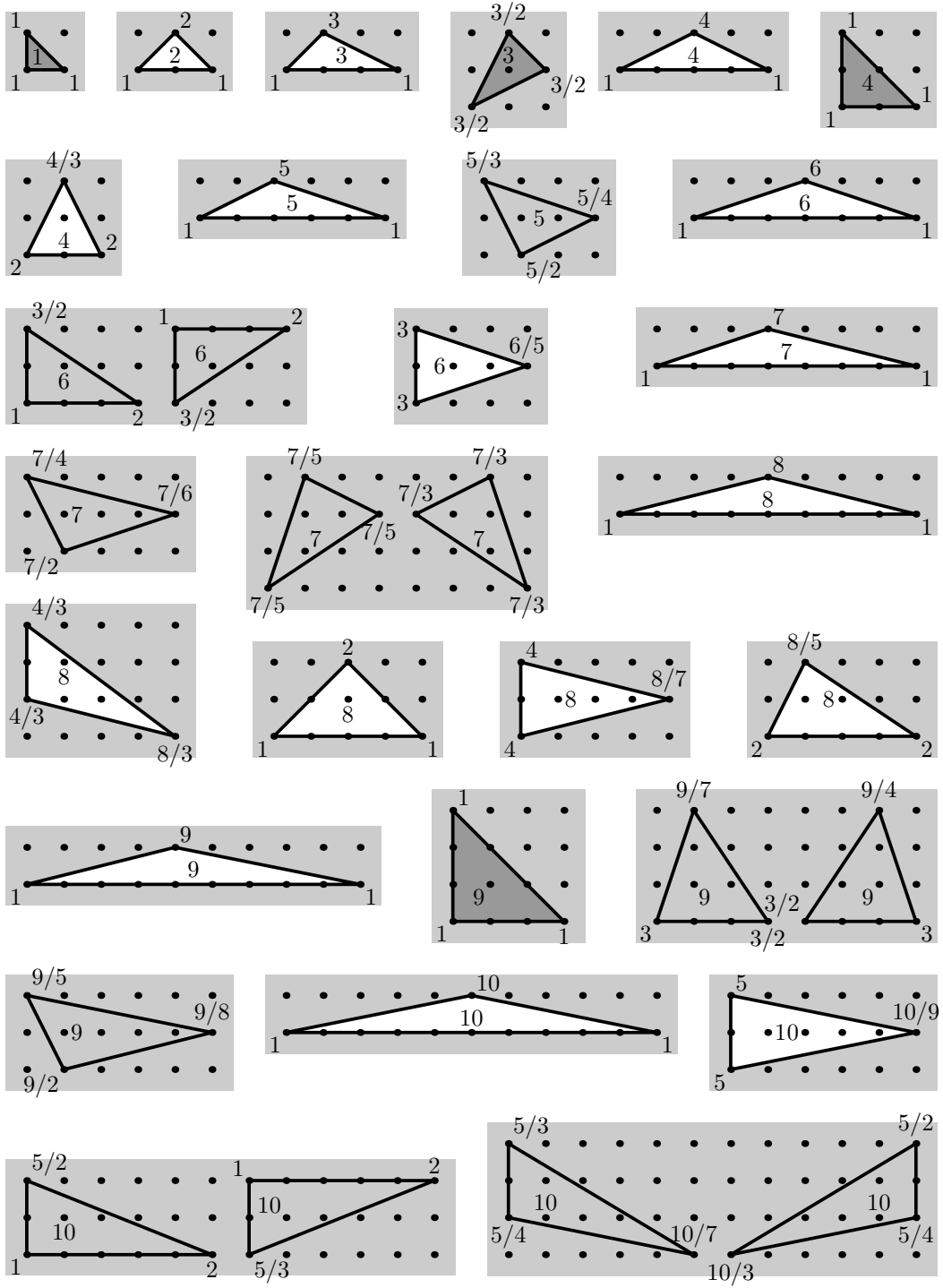}$$
\caption{List of lattice triangles of lattice volume less than or
equal 10.}\label{triang}
\end{figure}

Here we discuss the lattice-congruence criterions for lattice
triangles. By the first criterion of lattice-congruence for
lattice triangles we obtain that the number of all
lattice-congruence classes for lattice triangles with bounded
lattice area is finite. Further we introduce the complete list of
lattice-congruence classes for triangles with lattice area less
then or equal to 10. There are exactly $33$ corresponding
lattice-congruence classes, shown on Figure~\ref{triang}.

\vspace{2mm}

{\bf On criterions of lattice triangle lattice-congruence.} We
start with the study of lattice analogs for the first, the second,
and the third Euclidean criterions of triangle congruence.

\begin{statement}{\bf(The first criterion of lattice triangle lattice-congruence.)}
Consider two lattice triangles $\triangle ABC$ and $\triangle
A'B'C'$. Suppose that the edge $AB$ is lattice-congruent to the
edge $A'B'$, the edge $AC$ is lattice-congruent to the edge
$A'C'$, and the ordinary angle $\angle CAB$ is lattice-congruent
to the ordinary angle $\angle C'A'B'$, then the triangle
$\triangle A'B'C'$ is lattice-congruent to the triangle
$\triangle ABC$. \qed
\end{statement}

It turns out that the second and the third criterions taken from
Euclidean geometry do not hold. The following two examples
illustrate these phenomena.

\begin{example}
The second criterion of triangle lattice-congruence does not hold
in lattice geometry. On Figure~\ref{examples1} we show two
lattice triangles $\triangle ABC$ and $\triangle A'B'C'$. The
edge $AB$ is lattice-congruent to the edge $A'B'$ (here $\il
(A'B')=\il(AB)=4$). The ordinary angle $\angle ABC$ is
lattice-congruent to the ordinary angle $\angle A'B'C'$ (since
$\angle ABC \cong \angle A'B'C'\cong \iarctan (1)$), and  the
ordinary angle $\angle CAB$ is lattice-congruent to the ordinary
angle $\angle C'A'B'$ (since $\angle CAB \cong \angle C'A'B'\cong
\iarctan (1)$), The triangle $\triangle A'B'C'$ is not
lattice-congruent to the triangle $\triangle ABC$, since $\is
(\triangle ABC) =4$ and $\is (\triangle A'B'C') =8$.
\begin{figure}[h]
$$\epsfbox{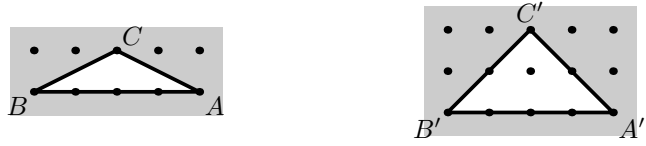}$$
\caption{The second criterion of triangle lattice-congruence does
not hold.}\label{examples1}
\end{figure}
\end{example}

\begin{example}
The third criterion of triangle lattice-congruence does not hold
in lattice geometry. On Figure~\ref{examples2} we show two
lattice triangles $\triangle ABC$ and $\triangle A'B'C'$. All
edges of both triangles are lattice-congruent (of length one),
but the triangles are not lattice-congruent, since $\is
(\triangle ABC) =1$ and $\is (\triangle A'B'C') =3$.
\begin{figure}[h]
$$\epsfbox{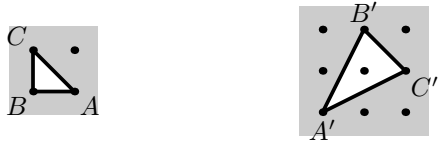}$$
\caption{The third criterion of triangle lattice-congruence does
not hold.}\label{examples2}
\end{figure}
\end{example}

Instead of the second and the third criterions there exists the
following additional criterion of lattice triangles
lattice-congruence.

\begin{statement}{\bf(An additional criterion of lattice triangle in\-te\-ger-con\-gru\-ence.)}
Consider two lattice triangles $\triangle ABC$ and $\triangle
A'B'C'$ of the same lattice area. Suppose that the ordinary angle
$\angle ABC$ is lattice-congruent to the ordinary angle $\angle
A'B'C'$, the ordinary angle $\angle CAB$ is lattice-congruent to
the ordinary angle $\angle C'A'B'$, the ordinary angle $\angle
BCA$ is lattice-congruent to the ordinary angle $\angle B'C'A'$,
then the triangle $\triangle A'B'C'$ is lattice-congruent to the
triangle $\triangle ABC$. \qed
\end{statement}

In the following example we show that the additional criterion of
lattice triangle lattice-congruence is not improvable.

\begin{example}
On Figure~\ref{examples3} we show an example of two lattice
non-equivalent triangles $\triangle ABC$ and $\triangle A'B'C'$
of the same lattice area equals 4 and the same ordinary lattice
angles $\angle ABC$, $\angle CAB$, and $\angle A'B'C'$, $\angle
C'A'B'$ all lattice-equivalent to the angle $\iarctan(1)$, but
$\triangle ABC\not\cong \triangle A'B'C'$.
\begin{figure}[h]
$$\epsfbox{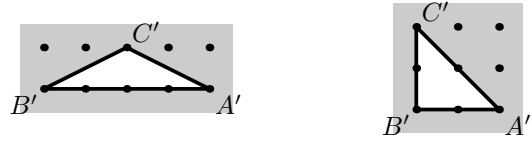}$$
\caption{The additional criterion of lattice triangle
lattice-congruence is not improvable.}\label{examples3}
\end{figure}
\end{example}

\vspace{2mm} {\bf Examples of lattice triangles.} First, we
define different types of lattice triangles.
\begin{definition}
The lattice triangle $\triangle ACB$ is called {\it dual} to the triangle $\triangle ABC$.\\
The lattice triangle is said to be {\it self-dual} if it is lattice-congruent to the dual triangle.\\
The lattice triangle is said to be {\it pseudo-isosceles} if it has at least two lattice-congruent angles.\\
The lattice triangle is said to be {\it lattice isosceles} if it is pseudo-isosceles and self-dual.\\
The lattice triangle is said to be {\it pseudo-regular} if all its ordinary angles are lattice-congruent.\\
The lattice triangle is said to be {\it lattice regular} if it is
pseudo-regular and self-dual.
\end{definition}

On Figure~\ref{triang} we show the complete list of $33$ triangles
representing all lattice-congruence classes of lattice triangles
with small lattice areas not greater than 10. We enumerate the
vertices of the triangle in the clockwise way. Near each vertex
of any triangle we write the tangent of the corresponding
ordinary angle. Inside any triangle we write its area. We draw
dual triangles on the same light gray area (if they are not
self-dual). Lattice regular triangles are colored in dark grey,
lattice isosceles but not lattice regular triangles are white,
and the others are light grey.

\section{Some unsolved question on lattice trigonometry.}

We conclude this paper with a small collection of unsolved
questions.

Let us start with some questions on elementary definitions of
lattice trigonometry. In this paper we do not show any
geometrical meaning of lattice cosine. Here arise the following
question.

\begin{problem}
Find a natural description of lattice cosine for ordinary lattice
angles in terms of lattice invariants of the corresponding
sublattices.
\end{problem}

This problem seems to be close to the following one.

\begin{problem}
Does there exist a lattice analog of the cosine formula for the
angles of triangles in Euclidean geometry?
\end{problem}

Now we formulate some problems on definitions of lattice
trigonometric functions for lattice irrational angles.

\begin{problem}
{\bf a).} Find a natural definition of lattice tangents for
lattice L-irrational angles, and
lattice LR-irrational angles.\\
{\bf b).} Find a natural definition of lattice sines and cosines
for lattice irrational angles?
\end{problem}

Let us continue with questions on lattice analogs of classical
trigonometric formulas for trigonometric functions of angles of
triangles in Euclidean geometry.

\begin{problem}
{\bf a).} Knowing the lattice trigonometric functions for lattice
angles $\alpha$, $\beta$ and integer $n$, find the explicit
formula for the lattice trigonometric functions
of the expanded lattice angle $\bar\alpha {+_n} \bar\beta$.\\
{\bf b).} Knowing the lattice trigonometric functions for a
lattice angle $\alpha$, an integer $m$, and positive integer $m$,
find the explicit formula for the lattice trigonometric functions
of the expanded lattice angle
$$
\sum\limits_{M,i=1}^{l}\bar{\alpha},
$$
where $M=(m,\ldots,m)$ is an $n$-tuple.
\end{problem}

The next problem is about a generalization of the statement of
Theorem~\ref{sum}b to the case of $n$ ordinary angles, that is
important in toric geometry and theory of multidimensional
continued fractions.

\begin{problem}
Find a necessary and sufficient conditions for the existence of
an $n$-gon with the given ordered sequence of ordinary lattice
angles $(\alpha_1, \ldots, \alpha_n)$ and the consistent sequence
of lattice lengths of the edges $(l_1, \ldots ,l_n)$ in terms of
continued fractions for $n\ge 4$.
\end{problem}

In Section~\ref{iic}, in particular, we introduced the definition
of the sums of any expanded lattice L-irrational angle with any
expanded lattice R-irrational angle.

\begin{problem}
Does there exist a natural definition of the sums of\\
{\bf a)} any expanded lattice LR-irrational angle and any expanded lattice (irrational) angle;\\
{\bf b)} any expanded lattice R-irrational angle and any expanded lattice (irrational) angle;\\
{\bf c)} any expanded lattice (irrational) angle and any expanded
lattice L-irrational angle?
\end{problem}

We conclude this paper with the following problem being actual in
the study of expanded lattice irrational angles.
\begin{problem}
Find an effective algorithm to verify whether two given
almost-positive signed lattice length-sine sequences define
lattice-congruent expanded lattice irrational angles, or not.
\end{problem}


\begin{thebibliography}{99}

\bibitem{Arn2}
V.~I.~Arnold, {\it Continued fractions}, M.: Moscow Center of
Continuous Mathematical Education, 2002.
\bibitem{Arn5}
V.~I.~Arnold, {\it Statistics of integer convex polygons}, Func.
an. and appl., v.~14(1980), n.~2, pp.~1--3.
\bibitem{Bar2}
I.~B\'ar\'any, A.~M.~Vershik, {\it On the number of convex lattice
polytopes}, Geom. Funct. Anal. v.~2(4), 1992, pp.~381--393.
\bibitem{Dan}
V.~I.~Danilov, {\it The geometry of toric varieties}, Uspekhi
Mat. Nauk, v.~33(1978), n.~2, pp.~85--134.
\bibitem{Ewa}
G.~Ewald, {\it Combinatorial Convexity and Algebraic Geometry},
Grad. Texts in Math. v.~168, Springer-Verlag, New York, 1996.
\bibitem{Ful}
W.~Fulton, {\it Introduction to Toric Varieties}, Annals of
Mathematics Studies; Princeton University Press, v.~131(1993),
\bibitem{Hin}
A.~Ya.~Hinchin, {\it Continued fractions}, M.: FISMATGIS, 1961.
\bibitem{Hir}
F.~Hirzebruch, {\it \"Uber vierdimensionale Riemannsche Fl\"achen
Mehrdeutinger analystischer Funktionen von zwei komplexen
Ver\"anderlichen}, Math. Ann. v.~126(1953), pp.~1--22.
\bibitem{Jun}
H.~W.~E.~Jung, {\it Darstellung der Funktionen eines
algebraischen K\"orpsers zweier unabh\"angigen Ver\"anderlichen
$x,y$ in der Umgebung einer Stelle $x=a$, $y=b$}, J.~Reine Angew.
Math., v.~133(1908), pp.~289--314.
\bibitem{Kar1}
O.~Karpenkov, {\it On tori decompositions associated with
two-dimensional continued fractions of cubic irrationalities},
Func. an. and appl., v.~38(2004), no~2, pp.~28--37.
\bibitem{KarReg}
O.~Karpenkov,
{\it Classification of lattice-regular lattice convex polytopes},\\
http://arxiv.org/abs/math.CO/0602193.
\bibitem{Kh1}
A.~G.~Khovanskii, A.~Pukhlikov, {\it Finitely additive measures
of virtual polytopes}, Algebra and Analysis, v.~4(2), 1992,
pp.~161--185; translation in  St.Petersburg Math. J., v.~4(2),
1993, pp.~337--356.
\bibitem{Kh2}
A.~G.~Khovanskii, A.~Pukhlikov, {\it A Riemann-Roch theorem for
integrals and sums of quasipolynomials over virtual polytopes},
Algebra and Analysis, v.~4(4), 1992, pp.~188--216; translation in
St. Petersburg Math. J., v.~4(4), 1993, pp.~789--812.
\bibitem{Kle1}
F.~Klein, {\it Ueber einegeometrische Auffassung der gew\"ohnliche
Kettenbruchentwicklung}, Nachr. Ges. Wiss. G\"ottingen Math-Phys.
Kl., 3(1891), 357--359.
\bibitem{Kon}
M.~L.~Kontsevich and Yu.~M.~Suhov, {\it Statistics of Klein
Polyhedra and Multidimensional Continued Fractions}, Amer. Math.
Soc. Transl., v.~197(2), (1999) pp.~9--27.
\bibitem{Kor2}
E.~I.~Korkina, {\it Two-dimensional continued fractions. The
simplest examples}, Proceedings of V.~A.~Steklov Math. Ins., v.
209(1995), pp. 243--166.
\bibitem{Lac2}
G.~Lachaud, {\it Voiles et Poly\`edres de Klein}, preprint n
95--22, Laboratoire de Math\'e\-matiques Discr\`etes du C.N.R.S.,
Luminy (1995).
\bibitem{Oda}
T.~Oda, {\it Convex bodies and Algebraic Geometry, An
Introduction to the Theory of Toric Varieties}, Springer-Verlag,
Survey in Mathemayics, 15(1988).
\bibitem{P-P}
P.~Popescu-Pampu, {\it The Geometry of Continued Fractions and
The Topology of Surface Singularities}, ArXiv:math.GT/0506432,
v1(2005).
\end{thebibliography}
\end{document}